\documentclass[r]{siamart171218}
\usepackage[english]{babel}
\usepackage{amsmath}
\usepackage{amssymb}
\usepackage{amsfonts}
\usepackage{array}
\usepackage{graphicx}
\usepackage{epsfig}
\usepackage{float}
\usepackage{color}
\usepackage{enumitem}  
\usepackage{epstopdf}
\usepackage{stmaryrd}
\usepackage{standalone}
\usepackage{tikz, subfigure}
\usepackage{yfonts}

\usetikzlibrary{shapes, patterns, arrows}
\usetikzlibrary{calc}
\usepackage{pgfplots}
\usepgfplotslibrary{fillbetween}

\usepackage{booktabs}  
\usepackage{capt-of}
\usepackage{multirow}

\newcommand{\kent}[1]{{\color{black}#1}}
\newcommand{\miro}[1]{{\color{black}#1}}
\newcommand{\paolo}[1]{{\color{black}#1}}
\title{Analysis and approximation of mixed-dimensional PDEs on 3D-1D domains coupled with Lagrange multipliers}

\author{M. Kuchta\thanks{Simula Research Laboratory, Oslo, Norway (\email{miroslav@simula.no})}\and F. Laurino\thanks{MOX, Department of Mathematics, Politecnico di Milano, Italy (\email{federica.laurino@polimi.it})} \and K.A. Mardal\footnotemark[1]$^{\;,}$\thanks{University of Oslo, Department of Mathematics (\email{kent-and@math.uio.no})} \and P. Zunino\thanks{MOX, Department of Mathematics, Politecnico di Milano, Italy (\email{paolo.zunino@polimi.it}), Corresponding Author}}

\begin{document}

\maketitle

\begin{abstract}
Coupled partial differential equations defined on domains with different dimensionality are usually called mixed dimensional PDEs. We address mixed dimensional PDEs on three-dimensional (3D) and one-dimensional domains, giving rise to a 3D-1D coupled problem.
Such problem poses several challenges from the standpoint of existence of solutions and numerical approximation. For the coupling conditions across dimensions, we consider the combination of essential and natural conditions, basically the combination of Dirichlet and Neumann conditions. To ensure a meaningful formulation of such conditions, we use the Lagrange multiplier method, suitably adapted to the mixed dimensional case. The well posedness of the resulting saddle point problem is analyzed. Then, we address the numerical approximation of the problem in the framework of the finite element method. The discretization of the Lagrange multiplier space is the main challenge. Several options are proposed, analyzed and compared, with the purpose to determine a good balance between the mathematical properties of the discrete problem and flexibility of implementation of the numerical scheme. The results are supported by evidence based on numerical experiments.
\end{abstract}



 

\def\ud{u_{\odot}}
\def\vd{v_{\odot}}
\def\ld{\lambda_{\odot}}
\def\udh{u_{\odot h}}
\def\udfh{u_{\odot \fkh}}
\def\vdh{v_{\odot h}}
\def\vdfh{v_{\odot \fkh}}
\def\ldh{\lambda_{\odot h}}
\def\md{\mu_{\odot}}
\def\mdh{\mu_{\odot h}}
\def\lld{l_{\odot}}
\def\uf{u_{\ominus}}
\def\up{u_{\oplus}}
\def\eps{\epsilon}
\def\nn{\boldsymbol n}
\def\rr{\boldsymbol r}
\def\RR{\boldsymbol R}
\def\kk{\boldsymbol k}
\def\ss{\boldsymbol s}
\def\uu{\boldsymbol u}
\def\vv{\boldsymbol v}
\def\xx{\boldsymbol x}
\def\bu{\overline{u}}
\def\bv{\overline{v}}
\def\tu{\widetilde{u}}
\def\tv{\widetilde{v}}
\def\TT{\boldsymbol T}
\def\NN{\boldsymbol N}
\def\BB{\boldsymbol B}
\def\ttu{\widetilde{\widetilde{u}}}
\def\ttv{\widetilde{\widetilde{v}}}
\def\cv{\check{v}}
\def\mesh{{\cal T}^h}
\def\ball{{\cal B}}
\def\R{\mathbb{R}}
\def\D{\mathcal{D}}
\def\DD{\partial\mathcal{D}}
\def\trace{{\mathcal{T}_\Gamma}}
\def\mtrace{{\overline{\mathcal{T}}_\Lambda}}
\def\ext{\mathcal{E}_\Gamma}
\def\ide{\mathcal{I}}
\def\ii{\hat{\imath}}
\def\fkh{\mathfrak{h}}
\newcommand{\avrd}[1]{\overline{\overline{#1}}}
\newcommand{\avrc}[1]{\overline{#1}}
\newcommand{\refe}[1]{{#1}_{\mathrm{ref}}}
\newcommand{\norm}[1]{\lVert{#1}\rVert}

\newcommand{\vertiii}[1]{{\left\vert\kern-0.25ex\left\vert\kern-0.25ex\left\vert #1 
    \right\vert\kern-0.25ex\right\vert\kern-0.25ex\right\vert}}

\newtheorem{thm}{Theorem}[section]
\newtheorem{prop}{Property}[section]
\theoremstyle{remark}
\newtheorem{remark}{Remark}[section]
 
\section{Introduction}\label{sec:intro}

In this study we consider coupled partial differential equations on domains with mixed dimensionality, in particular we address the  3D-1D case. The mathematical structure of such problems can be represented by the following formal equations:
\begin{subequations}
\label{eq:pde}
\begin{align}
\label{eq:pde1}
  -\Delta u + u + \lambda \delta_{\Lambda} &= f &\mbox{ in } \Omega,\\
\label{eq:pde2}
 d_s^2 \ud + \ud - \lambda &= g &\mbox{ on } \Lambda,\\
\label{eq:coupling}
\mathcal{T}_\Lambda u - \ud  &=  \paolo{q} &\mbox{ on } \Lambda . 
\end{align}
\end{subequations}
Problem \eqref{eq:pde} can be described as an example of \emph{mixed dimensional PDEs}.
Here, $u$, $\ud$,  $\lambda$ are unknowns,  $\Omega$ is a bounded domain in $\R^3$, whereas $\Lambda \subset \Omega$ is a 1D structure
parameterized in terms of $s$ and $d_s$ is the derivative with respect to $s$. 
The term $\lambda\delta_{\Lambda}$ is a Dirac measure such that 
$\int_{\Omega}\lambda(x)\delta_{\Lambda}v(x)\,\mathrm{d}x=\int_{\Lambda}\lambda(t)v(t) \,\mathrm{d}t$
for a continuous function $v$ and $\mathcal{T}_\Lambda: \Omega\rightarrow\Lambda$ is a suitable restriction operator from 3D to 1D. 

Using models based on mixed dimensional PDEs is motivated by the fact that many problems in geo- and biophysics are characterized by slender cylindrical structures coupled to a larger 3D body, where the characteristic transverse length scale of the slender structure is many orders of magnitude smaller than the longitudinal length.  For example, in geophysical applications the radii of wells are often of the order of 10 cm while the length may be several kilometers~\cite{Peaceman1978183, Peaceman1983531}. Similarly, in applications involving the blood flow and oxygen transport of the micro-circulation the capillary radius is a few microns, while simulations are often performed on mm to cm scale, with thousands of vessels~\cite{berg2020modelling,fang2008oxygen,gould2017capillary, secomb2004green}. Finally, in neuro-science applications a neuron has width of a few microns, while  its length is much longer. For example, an axon of a motor neurons may be as long as a meter.  Hence,  at least 4 orders of magnitude in difference in transverse and longitudinal direction is common in both geo-physics, bio-mechanics and neuro-science. Meshes dictated by resolving the transverse length scale in 3D would then possibly lead to the order of $10^{12}$ degrees of freedom. Even if adaptive and strongly anisotropic meshes are allowed for, the computations quickly become demanding if many slender structures and their interactions are under study.  

From a mathematical standpoint, the challenge involved in problem \eqref{eq:pde} is that neither $\mathcal{T}_\Lambda$ nor $\delta_\Lambda$ are well defined. That is, without extra regularity,  solutions of elliptic PDEs only have well defined traces of co-dimension one. Here, $\mathcal{T}_\Lambda$ is of co-dimension two, mapping functions defined on a domain in 3D to functions defined along a 1D curve. 
The challenge of coupling PDEs on domains with high dimensionality gap has recently attracted the attention of many researchers. The sequence of works by D'Angelo, \cite{DAngelo,d2012finite,d2008coupling} have remedied the well-posedness by weakening the solution concept. The approach naturally leads to non-symmetric formulations. An alternative approach is to decompose the solution into smooth and non-smooth components, where the non-smooth component may be represented in terms of Green's functions, and then consider the well-posedness of the smooth component \cite{gjerde2019singularity}. The numerical approximation of such equations has been also studied in a series of works. The consistent derivation of numerical approximation schemes for PDEs in mixed dimension is addressed in \cite{Nordbotten}. Concerning approximability,  elliptic equations with Dirac sources represent an effective prototype case that has been addressed in \cite{Bertoluzza201897,koppl2016local,koppl2014optimal}, where the optimal a-priori error estimates for the finite element approximation are derived. Furthermore, the interplay between the mathematical structure of the problem and solvers, as well as preconditioners for its discretization has been studied in details  in \cite{kuchta2016preconditioners} for the solution of 1D differential equations embedded in 2D, and more recently extended to the 3D-1D case in \cite{KMM2}.

Stemming from this literature, in this work we adopt and analyze a different approach, closely related to \cite{KVWZ,laurino_m2an}. That is, we exploit the fact that $\Lambda$ is not strictly a 1D curve, but rather a very thin 3D structure with a cross-sectional area far below from what can be resolved. With this additional assumption, we show that robustness with respect to the cross-sectional area can be restored. The major novelty with respect to the previous works is that we address essential type coupling conditions, namely Dirichlet-Neumann conditions, rather than natural type ones, such as the Neumann-Robin or Robin-Robin cases. These coupling conditions pose additional difficulties as the conditions 
are not a natural part of the weak formulation of the problem. We overcome this difficulty by  resorting to a weak formulation of the Dirichlet-Neumann coupling conditions across dimensions by using Lagrange multipliers.

\paolo{Although the focus of the present work is mostly about analysis and approximation of the proposed approach, we stress that it aims to build the mathematical foundations to tackle various applications involving 3D-1D mixed dimensional PDEs, such as FSI of slender bodies \cite{Mori2019887}, microcirculation and lymphatics \cite{Possenti2019101,Vinje2019}, subsurface flow models with wells \cite{Cerroni2019} and the electrical activity of neurons.}

\section{Preliminaries}\label{sec:setting}

\begin{figure}[t!]
\begin{center}
\includegraphics[width=0.35\textwidth,angle=-90]{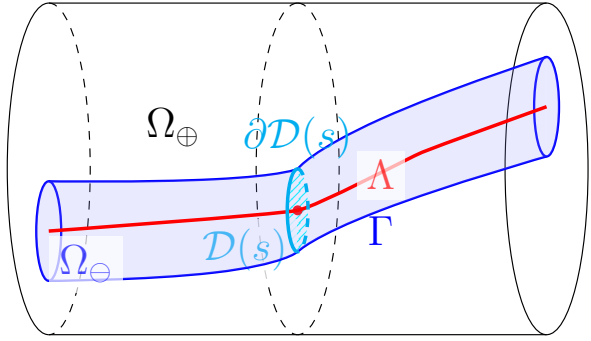}
\end{center}
\caption{Geometrical setting of the problem}
\label{fig1}
\end{figure}

Let the domain $\Omega\subset \R^3 $ be convex and composed of two parts, 
$\Omega_{\ominus}$ and $\Omega_{\oplus}:=\Omega\setminus\overline{\Omega}_{\ominus}$. 
Let 
$\Omega_{\ominus}$ be a \emph{generalized cylinder}, c.f.  \cite{MR1940257}, 
that is; the swept volume of a two dimensional set, $\partial\mathcal{D}$, moved along a curve, $\Lambda$, in the three-dimensional domain, $\Omega$, see for Figure \ref{fig1} for an illustration. 
In detail, the curve 
$\Lambda = \{\boldsymbol \lambda(s), \ s\in(0,S)\}$,  where 
$\boldsymbol \lambda(s) = [\xi(s), \tau(s), \zeta(s)], \ s\in(0,S)$ is a $\mathcal{C}^2$-regular curve in the three-dimensional domain $\Omega$.
For simplicity, let us assume that $\|\boldsymbol \lambda^\prime(s)\|=1$ such that the arc-length and the coordinate $s$ coincide.
Further, let $\mathcal{D}(s) = [x(r,t),y(r,t)]: (0,R(s))\times(0,T(s)) \rightarrow \R^2$ be a parametrization of the cross section
and $\Gamma$ be the lateral surface of $\Omega_{\ominus}$, i.e.
$\Gamma=\{ \ \partial \mathcal{D}(s) \  | \ s \in \Gamma\}$,
while the upper and lower faces of $\Omega_{\ominus}$ belong to $\partial\Omega$. 
We assume that $\Omega_{\ominus}$ crosses $\Omega$ from side to side.
 \paolo{Finally, $|\cdot|$ denotes the Lebesgue measure of a set, e.g. $|\mathcal{D}(s)|$ is the cross-sectional area of the cylinder. In general, $|\D(s)|$ must be strictly positive and bounded.}
According to the geometrical setting, we will denote with $v,\,v_\oplus,\,v_\ominus,\,\vd$,
functions defined on $\Omega,\,\Omega_{\oplus},\,\Omega_{\ominus},\,\Lambda$, respectively.

Let $D$ be a generic regular bounded domain in $\R^3$ and $X$ be a Hilbert space defined on $D$. Then  $(\cdot,\cdot)_X$ and $\|\cdot\|_X$ denote the inner product and norm of $X$, respectively.
The duality pairing between the $X$ and its dual $X^*$ is denoted as $\langle\cdot,\cdot\rangle$.
Let $(\cdot,\cdot)_{L^2(D)}$,  $(\cdot,\cdot)_D$ or simply $(\cdot,\cdot)$ be the $L^2(D)$ inner product on $D$.
We use the standard notation $H^q(D)$ to denote the Sobolev space of functions on $D$ with all derivatives up to the order $q$ in $L^2(D)$.
The corresponding norm is $\|\cdot\|_{H^q(D)}$ and the seminorm is $\vert\cdot\vert_{H^q(D)}$.
The space $H^q_0(D)$ represents the closure in $H^q(D)$ of smooth functions with compact support in $D$.

Let $\Sigma$ be a Lipschitz co-dimension one subset of $D$. We denote with $\mathcal{T}_\Sigma: H^q(D) \rightarrow H^{q-\frac 1 2}(\Sigma)$ the trace operator from \kent{$D$} to $\Sigma$.
The space of functions in $H^{\frac 1 2}(\Sigma)$ with continuous extension by zero outside $\Sigma$ is denoted 
$H^{\frac 1 2}_{00}(\Sigma)$ and we remark that 
$H^{\frac 1 2}_{00}(\Sigma) = \mathcal{T}_\Sigma H^1_0(D)$
and  $H^{-\frac 1 2}(\Sigma) = (H^{\frac 1 2}_{00}(\Sigma))^*$

\kent{
We will frequently use inner products and norms that are weighted. The
$L_2$ and $H^1$ inner products
weighted by a scalar function 
$w$, which is strictly positive and bounded almost everywhere,
are defined as 
follows
\begin{equation*}
(u,v)_{L^2(\Omega), w} = 
\int_{\Omega} w \, u \, v d\omega \ \mbox{   and   } \
(u,v)_{H^1(\Omega), w} = 
\int_{\Omega}  w \, u \, v d\omega
+ \int_{\Omega} w \, \nabla u \cdot \nabla v d\omega
\end{equation*}
whereas a weighted fractional
space \miro{$H^s(\Gamma; w)$} is defined in terms of the interpolation of the corresponding weighted spaces. 
For the norm of such spaces, we introduce 
the Riesz map $S$ such that 
for $u, v \in H^1(\Gamma)$ we have 
\[
(S u,v)_{H^1(\Gamma),w} = (u,v)_{L^2(\Gamma), w}.
\]
Then $S$ is a compact self-adjoint operator. 
Assuming that $\{\lambda_k\}_k$ is the set of eigenvalues, 
$\{\phi_k\}_k$ the set of eigenvectors of $S$ 
orthonormal with respect to the inner 
product $(\cdot, \cdot)_{L^2(\Gamma), w}$ and $u\in H^1(\Gamma)$ can be expressed 
as $u=\sum_k c_k \phi_k$ 
then
\[
\|u\|^2_{H^s(\Gamma), w} = \sum_k \lambda_k^{-s} c_k^2 .   
\]
The space $H_{00}^s(\Gamma; w)$ is defined analogously, but with $S$ above defined in terms of the $H^1_0$ inner product. Owing to the positivity and boundedness of $w$ the weighted spaces equal the corresponding non-weighted spaces as sets, but their norms are different.
 }

Central in our analysis are the transverse averages $\avrc{w},\,\avrd{w}$ defined as, 
\begin{gather*}
\avrc{w}(s) = |\partial\mathcal{D}(s)|^{-1} \int_{\partial\mathcal{D}(s)} w d\gamma \quad \mbox{and} \quad 
\avrd{w}(s) = |\mathcal{D}(s)|^{-1} \int_{\mathcal{D}(s)} w d\sigma, 
\end{gather*}
where $d\omega, \ d\sigma, \ d\gamma$ are the generic volume, surface and curvilinear Lebesgue measures.
Clearly, 
\begin{gather*}
\int_{\Omega_{\ominus}} w d\omega 
= \int_\Lambda \int_{\mathcal{D}(s)} w d\sigma ds
= \int_\Lambda |\mathcal{D}(s)|\avrd{w}(s) ds\, 
\\
\int_{\partial\Omega_{\ominus}} w d\sigma 
= \int_\Lambda \int_{\partial\mathcal{D}(s)} w d\gamma ds
= \int_\Lambda  |\partial \mathcal{D}(s)| \avrc{w}(s) ds\,. 
\end{gather*}
Analogously, for functions defined on $\Lambda$ and $\Omega_\ominus$ respectively, 
we let  $d_s$ and $\partial_s$ be the ordinary and partial derivative with respect to the arclength.

The operator obtained from a combination of the average operator $\avrc{(\cdot)}$ with the trace on $\Gamma$ will be denoted with $\mtrace = \avrc{(\cdot)} \circ \trace$, as it maps functions on $\Omega$ to functions on $\Lambda$.
Further, let the extension operator $\mathcal{E}_{\Gamma}: H^{\frac 1 2}_{00}(\Lambda) \rightarrow H^{\frac 1 2}_{00}(\Gamma)$ be defined such that $(\mathcal{E}_{\miro{\Gamma}}\vd)(x)=\vd(s)$, for any $x\in \DD(s)$. Then, the following identity shows that the transversal uniform extension operator is the inverse of the transversal average,
\begin{equation}\label{eq:ext}
\langle \mtrace u, \vd \rangle_{\Lambda, |\DD|} 
= \int _{\Lambda} |\DD| \left( \frac{1}{|\DD|} \int_{\DD} \trace u \, d\gamma \right) \vd \, ds 
= \langle \trace u, \ext \vd\rangle_{\Gamma}\,.
\end{equation}



\kent{With the above notation we are now able to formulate the precise weak formulations of the problems \eqref{eq:pde}, which we will call the \textbf{Problem 3D-1D-1D}. The problem} reads:
\paolo{given $f\in L^2(\Omega), \ g \in L^2(\Omega_{\ominus}), \ q \in H^\frac12_{00}(\Gamma)$}
find $u \in H^1_0(\Omega),\ \ud \in H^1_0(\Lambda), \ \ld \in H^{-\frac 1 2}(\Lambda)$, such that
\begin{subequations}\label{eq:problem2}
\begin{align}
\label{eq:problem2a}
(u,v)_{H^1(\Omega)} 
+  \langle \mtrace v, \ld \rangle_{\Lambda, |\DD|} &= 
(f,v)_{L^2(\Omega)} &\forall v \in H^1_0(\Omega)\, ,
\\
\label{eq:problem2b}
(\ud,\vd)_{H^1(\Lambda),|\D|} - \langle  \vd, \ld \rangle_{\Lambda, |\DD|} 
&=  (\avrd{g},\vd)_{L^2(\Lambda),|\D|}
&\forall \vd \in H^1_0(\Lambda)\, ,
\\
\label{eq:problem2c}
\langle \mtrace u -   \ud, \md \rangle_{\Lambda,|\DD|} 
&= \paolo{\langle \avrc{q},\md\rangle_{\Lambda,|\DD|}}
&\forall \md \in H^{-\frac 1 2}(\Lambda)\,.
\end{align}
\end{subequations}

In addition to the 3D-1D-1D problem we will also consider an 
intermediate problem where the 3D and 1D problems are coupled at an
intermediate 2D surface  encapsulating the 1D structure. The strong form is:  
\paolo{
\begin{subequations}
\label{inteq:pde}
\begin{align}
\label{inteq:pde1}
  -\Delta u + u + \lambda \delta_{\Gamma} &= f &\mbox{ in } \Omega,\\ 
\label{inteq:pde2}
 d_s^2 \ud + \ud -  \avrc{\lambda} &= \avrd{g} &\mbox{ on } \Lambda,\\
\label{inteq:coupling}
\mathcal{T}_\Gamma u - \mathcal{E}_\Gamma \ud  &=  q &\mbox{ on } \Gamma . 
\end{align}
\end{subequations}
}
The corresponding weak formulation of \eqref{inteq:pde}, referred to as the \textbf{Problem 3D-1D-2D}, reads: 
\paolo{given $f\in L^2(\Omega), \ g \in L^2(\Omega_{\ominus}), \ q \in H^\frac12_{00}(\Gamma)$}
find $u \in H^1_0(\Omega), \ \ud \in H^1_0(\Lambda), \ \lambda \in H^{-\frac 1 2}(\Gamma ) $ such that
\begin{subequations}\label{eq:problem1}
\begin{align}
(u,v)_{H^1(\Omega)} + \langle \trace v, \lambda \rangle_\Gamma &= (f,v)_{L^2(\Omega)} &\forall v \in H^1_0(\Omega)\, ,
\\
 (\ud,\vd)_{H^1(\Lambda),|\D|}  -  \langle \ext \vd, \lambda \rangle_\Gamma
&=  (\avrd{g},\vd)_{L^2(\Lambda),|\D|} &\forall\vd \in H^1(\Lambda)\, ,
\\
\langle \trace u - \ext \ud , \mu \rangle_\Gamma 
&= \paolo{\langle q,\md\rangle_\Gamma} 
&\forall \mu \in H^{-\frac 1 2}(\Gamma)\,.
\end{align}
\end{subequations}

We conclude this section with the analysis of a fundamental property for the problem formulation that we will address, namely, the characterization of the regularity of the operator $\mtrace$. More precisely we aim to show that $\mtrace: H^1_0(\Omega) \rightarrow H^\frac 1 2_{00}(\Lambda)$.
This is a consequence of the following lemma.

\begin{lemma}\label{lemma:H12norm}
Let $\Gamma$ be a tensor product domain, $\Gamma= (0,X) \times (0,Y)$. For any regular $u(x,y)$ in $\Gamma$, let $\avrc{u}(x)=\frac 1Y \int _0^Y u(x,y)\, dy$. Then, for any $u\in H_{00}^{\frac 12}(\Gamma)$, $\avrc{u}(x)\in H_{00}^{\frac 12}((0,X))$. 
Moreover, if $u(x,y)\in H^{\frac 12}_{00}(\Gamma)$ is constant with respect to $y$, namely $u(x,y)=u(x)$, then 
\begin{equation*}
\|u\|_{H^{\frac 12}_{00}(\Gamma)}=Y \|u\|_{H^{\frac 12}_{00}(0,X)}.
\end{equation*}
\end{lemma}
The proof of \ref{lemma:H12norm} is based on the representation of fractional norms in terms of the spectrum of the Laplace operator and subsequent standard arguments in harmonic analysis. The full proof is reported in the appendix for the sake of clarity.



Under the geometric assumptions stated above for $\Omega,\,\Gamma,\,\Lambda$, Lemma \ref{lemma:H12norm} implies the following result.
\begin{corollary}\label{corollary:H12norm}
If $u\in H^{\frac 12}_{00}(\Gamma)$ then $\avrc{u}\in H^{\frac 12}_{00}(\Lambda)$ and there exists a constant $C_\Gamma$, \paolo{bounded independently of $\D$ and $\DD$, such that}
\begin{equation*}
\|\avrc{u}\|_{H^{\frac 12}_{00}(\Lambda), |\DD|} \leq C_\Gamma \|u\|_{H^{\frac 12}_{00}(\Gamma)}.
\end{equation*}
\end{corollary}
\paolo{
Furthermore, from the
above Corollary, it is clear that
$\mtrace: H^1_0(\Omega) \rightarrow H^\frac 1 2_{00}(\Lambda)$.
}

\section{Saddle-point problem analysis}
Let $a: X \times X \rightarrow \mathbb{R}$ and $b: X\times Q \rightarrow \mathbb{R}$ be bilinear forms. Let us consider a general saddle point problem of the form: find $u\in X$, $\lambda\in Q$ s.t.
\begin{eqnarray}\label{eq:saddle-point}
\begin{aligned}
a(u,v)+b(v,\lambda)&=c(v), &\forall v\in X,\\
b(u,\mu)&=d(\mu), &\forall \mu\in Q.
\end{aligned}
\end{eqnarray}
The Brezzi conditions \cite{MR365287} ensure that the problem \eqref{eq:saddle-point} is well-posed.
For our purpose here, we use the following relaxed version of the Brezzi conditions:  
\begin{theorem}{}\label{th:bnb}
Problem \eqref{eq:saddle-point} is well posed if the following conditions are satisfied
\begin{align}\label{BNB1}
a(u,u) &\ge \alpha \|u\|^2_{X}, & u\in X,\\
\label{BNB2}
a(u,v) &\le C \|u\|_{X}\|v\|_{X}, & u,v \in X,\\
\label{BNB3}
b(u,\mu) &\ge D \|u\|_{X} \|\mu\|_{Q}, & u\in X, \mu\in Q,\\
\label{BNB4}
\sup_{v\in X} \frac{b(v,\mu)}{\|v\|_{X}} &\geq \beta\|\mu\|_{Q}, & \mu \in Q .
\end{align}
Here $\alpha$, $\beta$, $C$, and $D$ are positive numbers. 
\end{theorem} 
Here, the coercivity condition \eqref{BNB1} applies to $X$, which is a relaxation of Brezzi's original conditions.

\subsection{Problem 3D-1D-2D}
We aim to find $u \in H^1_0(\Omega),\ \ud \in H_0^1(\Lambda), \ \lambda \in H^{-\frac12}(\Gamma )$,
solutions of \eqref{eq:saddle-point}, where
\begin{align*}
a([u, \ud], [v, \vd])&= (u,v)_{H^1(\Omega)} + (\ud,\vd)_{H^1(\Lambda),|\D|},
\\
b([v, \vd], \mu)&= \langle \trace v - \ext \vd , \mu \rangle_\Gamma,
\\
c([v,\vd])&= (f,v)_{L^2(\Omega)} + (\avrd{g},\vd)_{L^2(\Lambda),|\D|},
\\
d(\mu)&=\paolo{\langle q,\mu\rangle_\Gamma}.
\end{align*}
We prove that the conditions of Theorem \ref{th:bnb} are fulfilled choosing 
$X=H^1_0(\Omega) \times H^1_0(\Lambda)$, $Q=H^{-\frac 12}(\Gamma)$, where $X$  is equipped with the norm $\vertiii{[u,\ud ]}^2=\|u\|^2_{H^1(\Omega)} + \|\ud\|^2_{H^1(\Lambda),|\D|}$.

\begin{lemma}\label{lemma:prob1_boundedness} 
The Problem 3D-1D-2D is well-posed.
\end{lemma}
\begin{proof}
We need to establish the four Brezzi conditions. 
The bilinear form $a(\cdot \ , \ \cdot)$ is clearly bounded and coercive since 
for $u=\ud$, $v=\vd$
\begin{equation*}
a([u, \ud], [v, \vd]) 
= (u,v)_{H^1(\Omega)} + (\ud, \vd)_{H^1(\Lambda), |\D|}
= \|u\|^2_{H^1(\Omega)} + \|\ud\|^2_{H^1(\Lambda),|\D|}.
\end{equation*}
Furthermore, the bilinear form $b(\cdot \ , \ \cdot)$ is bounded because
\begin{align*}
b([v, \vd], \mu)&= \langle \trace v - \ext \vd , \mu \rangle_\Gamma
\leq \|\trace v - \ext \vd\|_{H^{\frac 12}_{00}(\Gamma)}\|\mu\|_{H^{-\frac 12}(\Gamma)}\\
&\leq \left(\|\trace v\|_{H^{\frac 12}_{00}(\Gamma)} + \|\ext \vd\|_{H^{\frac 12}_{00}(\Gamma)}\right)\|\mu\|_{H^{-\frac 12}(\Gamma)} \\
&\leq \left(C_T \|v\|_{H^1(\Omega)} + \|\ext \vd\|_{H^1(\Gamma)}\right)\|\mu\|_{H^{-\frac 12}(\Gamma)}\\
&\leq \left(C_T \|v\|_{H^1(\Omega)} + \left(\frac{\max |\DD|}{\min |\D|}\right)^{\frac 12} \|\vd\|_{H^1(\Lambda),|\D|}\right)\|\mu\|_{H^{-\frac 12}(\Gamma)}\\
&\leq \left( C_T + \left(\frac{\max |\DD|}{\min |\D|}\right)^{\frac 12}\right) \vertiii{[v,\vd]}\|\mu\|_{H^{-\frac 12}(\Gamma)}.
\end{align*}

\paolo{
To show the inf-sup condition, we will employ  
a lifting operator,  
$\mathcal{H}_\Omega$,
from $H^{1/2}_{00}(\Gamma)$ to $H^1(\Omega)$.
In \cite{sauter1999extension} 
it is established that extension operators for domains having small geometric details (see also \cite{laurino_m2an} for a direct application to this case) there exists a lifting operator $\mathcal{H}_\Omega$ from $H^{1/2}_{00}(\Gamma)$ to $H^1(\Omega)$ such that $\mathcal{H}_\Omega \xi = v$ for any $\xi \in H^{1/2}_{00}(\Gamma)$ with $v\in H^1(\Omega)$. Further, for this operator   there exists 
$\| \mathcal{H}_\Omega \| \in \mathbb{R}$ such that
$\|v \|_{H^1(\Omega)}\leq \|\mathcal{H}_\Omega\| \|\xi \|_{H^{1/2}_{00}(\Gamma)}$
where $\| \mathcal{H}_\Omega \|$ is a constant independent the (minimal) radius of $\Gamma$.
}


The inf-sup inequality is fulfilled, that is; 
we choose $\vd \in H^1_0(\Lambda)$ such that $\ext\vd =0$. Therefore,
\begin{equation*}
\sup _{\substack{v\in H^1_0(\Omega),\\ \vd \in H^1_0(\Lambda)}} \frac{ \langle \trace v  - \ext \vd, \mu\rangle_\Gamma}{\vertiii{[v, \vd]}} 
\geq \sup _{v\in H^1_0(\Omega)} \frac{ \langle \trace v, \mu \rangle_\Gamma}{\|v\|_{H^1(\Omega)}}.
\end{equation*}
We notice that the trace operator is surjective from $H^1_0(\Omega)$ to $H^{\frac12}_{00}(\Gamma)$. Indeed, $\forall \xi \in H^{\frac 12}_{00}(\Gamma)$, 
we  can find $v=\mathcal{H}_\Omega \xi$. Using the stability of the harmonic extension we obtain
\begin{equation}\label{infsup_traceop}
\sup _{v\in H^1_0(\Omega)} \frac{ \langle \trace v, \mu \rangle_\Gamma}{\|v\|_{H^1(\Omega)}}
\geq  \sup _{\xi \in H^{\frac 12}_{00}(\Gamma )} \frac{ \langle \xi , \mu \rangle_\Gamma}{\|\mathcal{H}_\Omega\| \|\xi\|_{H^{\frac 12}_{00}(\Gamma)}}
= \|\mathcal{H}_\Omega\|^{-1} \|\mu\|_{H^{-\frac 12}(\Gamma)},
\end{equation}
where in the last inequality we exploited the fact that $H^{-\frac 12}(\Gamma)=(H^{\frac 12 }_{00}(\Gamma))^*$. 
\end{proof}

\subsection{Problem 3D-1D-1D}
We aim to find $u \in H^1_0(\Omega),\ \ud \in H^1_0(\Lambda), \ \ld \in H^{-\frac12}(\Lambda)$,
solution of \eqref{eq:saddle-point} with
\begin{align*}
a([u, \ud], [v, \vd])&= (u,v)_{H^1(\Omega)} + (\ud,\vd)_{H^1(\Lambda),|\D|},
\\
b([v, \vd], \md)&=  \langle  \mtrace v - \vd, \md \rangle_{\Lambda, |\DD|},
\\
c([v,\vd])&= (f,v)_{L^2(\Omega)} + (\avrd{g},\vd)_{L^2(\Lambda),|\D|},
\\
d(\md)&=\paolo{\langle \avrc{q},\md\rangle_{\Lambda,|\DD|}}.
\end{align*}

We prove that the hypotesis of  Theorem \ref{th:bnb} are fulfilled with the following spaces $X=H^1_0(\Omega) \times H^1_0(\Lambda)$, $Q=H^{-\frac 12}(\Lambda)$.
Let us consider $X$ equipped again with the norm $\vertiii{[\cdot,\cdot ]}$ and  
$Q$ equipped with the norm $\|\cdot \|_{H^{-\frac 12}(\Lambda),|\DD|}$.
Then, we have the following lemmas.
\begin{lemma}
The Problem 3D-1D-1D is well-posed.
\end{lemma}
\begin{proof}
Again, 
\begin{equation*}
a([u, \ud], [v, \vd]) = (u,v)_{H^1(\Omega)} + (\ud, \vd)_{H^1(\Lambda), |\D|}.
\end{equation*}
The bound on $b(\cdot \ , \ \cdot)$ is established as
\begin{multline*}
b([v, \vd], \md)=  \langle  \mtrace v - \vd, \md \rangle_{\Lambda, |\DD|} 
\leq \|\mtrace v - \vd\|_{H^{\frac 12 }_{00}(\Lambda), |\DD|}\|\md\|_{H^{-\frac 12}(\Lambda),|\DD|}\\
\leq \left(\|\mtrace v\|_{H^{\frac 12 }_{00}(\Lambda), |\DD|}+\|\vd\|_{H^{\frac 12 }_{00}(\Lambda), |\DD|}\right)\|\md\|_{H^{-\frac 12}(\Lambda),|\DD|}\\
\leq \left(C_\Gamma \|\trace v\|_{H^{\frac 12 }_{00}(\Gamma)}+ \|\vd\|_{H^1(\Lambda), |\DD|}\right)\|\md\|_{H^{-\frac 12}(\Lambda),|\DD|} \\
\leq \left(C_\Gamma C_T\|v\|_{H^1(\Omega)}+ \left(\frac{\max |\DD|}{\min |\D|}\right)^{\frac 12} \|\vd\|_{H^1(\Lambda), |\D|}\right) \|\md\|_{H^{-\frac 12}(\Lambda),|\DD|} \\
\leq  \left( C_\Gamma C_T + \left(\frac{\max |\DD|}{\min |\D|}\right)^{\frac 12} \right) \vertiii{[v,\vd]}\|\md\|_{H^{-\frac 12}(\Lambda),|\DD|}.
\end{multline*} 
The inf-sup condition holds. 
We choose $\vd=0$ and obtain
\begin{equation*}
\sup _{\substack{v\in H^1_0(\Omega),\\ \vd \in H^1_0(\Lambda)}} \frac{ \langle \mtrace v - \vd, \md \rangle_{\Lambda,|\DD|}}{\vertiii{[v,\vd]}}
\geq \sup _{v\in H^1_0(\Omega)} \frac{ \langle \mtrace v, \md \rangle_{\Lambda,|\DD|}}{\|v\|_{H^1(\Omega)}}. 
\end{equation*}

For any $q \in H^{\frac 12}_{00}(\Lambda)$, we consider its uniform extension to $\Gamma$ named as $\ext q$
and then we consider the harmonic extension $v=\mathcal{H}_\Omega \ext q\in H^1_0(\Omega)$. It follows that $\mtrace v=q$. Therefore, 
\begin{equation*}
\sup _{v\in H^1_0(\Omega)}  \langle \mtrace v, \md \rangle_{\Lambda,|\DD|} \geq\sup_{q \in H^{\frac 12}_{00}(\Lambda)} \langle q, \md  \rangle_{\Lambda,|\DD|}\,.
\end{equation*}
Moreover, using Lemma \ref{lemma:H12norm} we obtain
\begin{equation*}
\|v\|_{H^1_0(\Omega)}\leq \|\mathcal{H}_\Omega\| \|\ext q\|_{H^{\frac 12}_{00}(\Gamma)}  = \|\mathcal{H}_\Omega\| \|q\|_{H^{\frac 12}_{00}(\Lambda),|\DD|}.
\end{equation*}
 Therefore, we conclude the proof with the following inequalities,
\begin{multline*}
\sup _{v\in H^1_0(\Omega)} \frac{ \langle \mtrace v, \md \rangle_{\Lambda,|\DD|}}{\|v\|_{H^1(\Omega)}}
\geq \sup _{q\in H^{\frac 12}_{00}(\Lambda)} \frac{ \langle q, \md \rangle_{\Lambda,|\DD|}}{\|v\|_{H^1(\Omega)}}
\\
\geq \frac{1}{\|\mathcal{H}_\Omega\|} \sup _{q\in H^{\frac 12}_{00}(\Lambda)} \frac{ \langle q, \md \rangle_{\Lambda,|\DD|}}{\|q\|_{H^{\frac 12}_{00}(\Lambda),|\DD|}} 
= \frac{1}{\|\mathcal{H}_\Omega\|} \|\md\|_{H^{-\frac 12}(\Lambda),|\DD|}.
\end{multline*}
\end{proof}


\section{Finite element approximation}
In this section we consider the discretization of the Problems 3D-1D-2D and 3D-1D-1D by means of the finite element method. We address two main challenges; first we aim to identify a suitable approximation space for the Lagrange multiplier and to analyze the stability of the discrete saddle point problem; second we aim to derive a stable discretization method that uses independent computational meshes for $\Omega$ and $\Lambda$, not necessarily conforming to $\Gamma$. Let us introduce a shape-regular triangulation $\mathcal{T}^{\Omega}_h$ of $\Omega$ and an admissible partition $\mathcal{T}^{\Lambda}_{\fkh}$ of $\Lambda$.
We analyze two different cases: the conforming case, where compatibility constraints are satisfied by $\mathcal{T}^{\Omega}_h$ and $\mathcal{T}^{\Lambda}_{h}$ with respect to $\Gamma$ and consequently $h=\fkh$; and the non conforming case, where it is possible to choose $\mathcal{T}^{\Omega}_h$ and $\mathcal{T}^{\Lambda}_{\fkh}$ arbitrarily.

The discrete equivalent of \eqref{eq:saddle-point} reads as finding $u_h\in X_h\subset X$, $\lambda_h\in Q_h\subset Q$ s.t.
\begin{eqnarray}\label{eq:saddle-point_discrete}
\begin{aligned}
a(u_h,v_h)+b(v_h,\lambda_h)&=c(v_h) &&\forall v_h\in X_h,\\
b(u_h,\mu_h)&=d(\mu_h) &&\forall \mu_h\in Q_h,
\end{aligned}
\end{eqnarray}
where with little abuse of notation we use $h$ as the sub-index for all the discretization spaces.
This discrete problem is well-posed if the \eqref{BNB1}-\eqref{BNB4} conditions applies to $X_h$ and $Q_h$. Since $X_h\subset X$ and $Q_h\subset Q$,  
\eqref{BNB1}-\eqref{BNB3} follow immediately and only the inf-sup condition needs consideration.  We summarize this proposition in the Corollary below. 
\begin{corollary}{\cite[Theorem 2.42]{MR2050138}}\\
Let
$X_h \subset X$, $Q_h \subset Q$, $a(\cdot, \cdot)$ and $b(\cdot, \cdot)$ satisfy the conditions  \eqref{BNB1}-\eqref{BNB3} then  
the problem \eqref{eq:saddle-point_discrete} is well-posed if 
the discrete counterpart of \eqref{BNB4} is satisfied, i.e.
there exists a constant $\beta_h>0$ such that
\begin{align}
\label{eq:infsup_discrete}
 \sup_{v_h\in X_h} \frac{b(v_h,\mu_h)}{\|v_h\|_{X}}\geq \beta_h \|\mu_h\|_{Q}\,, \kent{\quad \forall \mu_h \in Q_h.}
\end{align}
\end{corollary}

\subsection{Analysis of the case where $\mathcal{T}^{\Omega}_h$ conforms to $\mathcal{T}^{\Lambda}_\fkh$ and to $\Gamma$}
As conformity conditions between $\mathcal{T}^{\Omega}_h$, $\mathcal{T}^{\Lambda}_\fkh$ and $\Gamma$, we require that the intersection of $\mathcal{T}^{\Omega}_h$ and $\Gamma$ is made of entire faces of elements $K \in \mathcal{T}^{\Omega}_h$. Furthermore, we also set a restriction between $\mathcal{T}^{\Omega}_h$ and $\mathcal{T}^{\Lambda}_\fkh$. We assume that $\Lambda$ is a piecewise linear manifold. We want that the intersection of $\Gamma$ with any orthogonal plane to $\Lambda$ that crosses $\Lambda$ at the internal nodes of $\mathcal{T}^{\Lambda}_\fkh$, consists of entire edges of $\mathcal{T}^{\Omega}_h$. As a result of the latter condition we have $h=\fkh$. 
For this reason, we denote as $\mathcal{T}^{\Lambda}_h$ the mesh on $\Lambda$
from now on throughout this section.
 
\subsubsection{Problem 3D-1D-2D}
We denote by $X_{h,0}^k(\Omega)\subset H^1_0(\Omega)$, with $k>0$, the conforming finite element space of continuous piecewise polynomials of degree $k$ defined on $\Omega$ satisfying homogeneous Dirichlet conditions on the boundary and by $X_{h,0}^k(\Lambda)\subset H^1_0(\Lambda)$ the space of continuous piecewise polynomials of degree $k$ defined on $\Lambda$, satisfying homogeneous Dirichlet conditions on $\Lambda \cap \partial \Omega$. The space $Q_h$ must be suitably chosen such that \eqref{eq:infsup_discrete} holds. Let $Q_h$ be the trace space of $X_{h,0}^k(\Omega)$, namely the space of continuous piecewise polynomials of degree $k$ defined on $\Gamma$ which satisfy homogeneous Dirichlet conditions on $\partial \Omega$. As a result, $Q_h=X_{h,0}^k(\Gamma) \subset H^\frac12_{00}(\Gamma)$. The discrete version of the 3D-1D-2D problem is:
find $u_h \in X_{h,0}^k (\Omega) ,\, {\ud}_h \in X_{h,0}^k(\Lambda) ,\, \lambda_h \in Q_h \subset H^{-\frac12}(\Gamma )$, such that
\begin{subequations}
\begin{align}
&(u_h,v_h)_{H^1(\Omega)} + ({\ud}_h,{\vd}_h)_{H^1(\Lambda),|\D|}
+ \langle \trace v_h  - \mathcal{E}_{\Lambda} {\vd}_h, \lambda_h \rangle_\Gamma
\\
\nonumber
&\quad = (f,v_h)_{L^2(\Omega)} + (\avrd{g},{\vd}_h)_{L^2(\Lambda),|\D|} 
\quad \forall v_h \in X_{h,0}^k(\Omega), \ {\vd}_h \in X_{h,0}^k(\Lambda)\,,
\\
&\langle \trace u_h - \mathcal{E}_{\Lambda} {\ud}_h , \mu_h \rangle_\Gamma 
= \paolo{\langle q,\mu_h \rangle_\Gamma}
\quad \forall \mu_h \in Q_h\,.
\end{align}
\end{subequations}
In what follows, we analyze the well-posedness of the discrete problem. 
From now on, $C$ denotes a generic constant independent of the mesh size.

\begin{lemma}\label{lemma:prob1_orthproj}
Let $P_h: H^{\frac 12}_{00}(\Gamma) \rightarrow Q_h$ be the orthogonal projection operator defined  for any $v \in H^{\frac 12}_{00}(\Gamma)$ by
$(P_h v , \psi_h)_\Gamma= (v, \psi_h)_\Gamma$ for any $\psi_h \in Q_h$.
Then, $P_h$ is continuous on $H^{\frac 12}_{00}(\Gamma)$, namely
$\|P_h v\|_{H^{\frac 12}_{00}(\Gamma)} \leq C \|v\|_{H^{\frac 12}_{00}(\Gamma)}$.
\end{lemma}

\begin{proof}
We show that $P_h$ is continuous on $L^2(\Gamma)$ and on $H^1_0(\Gamma)$ following \cite[Section 1.6.3]{MR2050138}.  Then, Lemma \ref{lemma:prob1_orthproj} can be proved by interpolation between spaces, since $H^\frac12_{00}(\Gamma)$ can be seen as the interpolation space between $L^2(\Gamma)$ and $H^1_0(\Gamma)$. For the $L^2$-continuity, we exploit the fact that, from the definition of $P_h$, $(v-P_h v,P_h v)_{\Gamma}=0$.
Therefore, by Pythagoras identity,
\begin{equation*}
\|v\|^2_{L^2(\Gamma)} = \|v-P_h v\|_{L^2(\Gamma)}^2 + \|P_h v\|_{L^2(\Gamma)}^2 \geq \|P_h v\|^2 _{L^2(\Gamma)}.
\end{equation*}
Let us now consider $v\in H^1_0(\Gamma)$. The Scott-Zhang interpolation operator $SZ_h$ from $H^1_0(\Gamma)$ to $Q_h$ satisfies the following inequalities,
\begin{align}
\label{SZ_stability}
\|SZ_h v\|_{H^1(\Gamma)} &\leq C_1 \|v\|_{H^1(\Gamma)},
\\
\label{SZ_approx}
\|v -SZ_h v \|_{L^2(\Gamma)} &\leq C_2 h \|v\|_{H^1(\Gamma)}.
\end{align}
Therefore, using \eqref{SZ_stability}, \eqref{SZ_approx}, 
the $L^2$ stability of $P_h$ and the inverse inequality, we obtain,
\begin{multline*}
\|\nabla P_h v\|_{L^2(\Gamma)} 
\leq \|\nabla (P_h v - SZ_h v)\|_{L^2(\Gamma)} + \|\nabla SZ_h v\|_{L^2(\Gamma)}
\\
\leq  \|\nabla (P_h v - SZ_h v)\|_{L^2(\Gamma)} + C_1\|v\|_{H^1(\Gamma)}
= \frac{C_3}{h} \|P_h (v - SZ_h v)\|_{L^2(\Gamma)} + C_1\|v\|_{H^1(\Gamma)}
\\
\leq \frac{C_3}{h} \|v - SZ_h v\|_{L^2(\Gamma)} + C_1\|v\|_{H^1(\Gamma)}
\leq (C_2 C_3 +C_1) \|v\|_{H^1(\Gamma)},
\end{multline*}
from which we obtain the continuity in $H^1_0(\Gamma)$.
\end{proof}

\begin{lemma}\label{lemma:trspace_infsup} 
There exists a constant $\gamma >0$ such that for any $\mu_h\in Q_h$
\begin{equation*}
\sup_{\substack{q_h \in Q_h}} \frac{\langle q_h , \mu_h \rangle}{ \|q_h\|_{H^{\frac 12}_{00}(\Gamma)}} \geq \gamma \|\mu_h\|_{H^{-\frac 12}(\Gamma)}.
\end{equation*} 
\end{lemma}

\begin{proof}
From the continuous case, in particular from \eqref{infsup_traceop}, we have
\begin{equation*}
\|\mathcal{H}_\Omega\|^{-1} \|\mu_h\|_{H^{-\frac 12}(\Gamma)} \leq \sup_{\substack{v \in H^1_0(\Omega)}} \frac{\langle \trace v , \mu_h \rangle}{\|v\|_{H^1(\Omega)}} 
\quad \forall \mu_h \in Q_h,
\end{equation*}
and by the trace inequality $\|\trace v\|_{H^\frac 12 (\Gamma)} \leq C_T \|v\|_{H^1(\Omega)}$ (see \cite[7.56]{adams1975pure}), we obtain 
\begin{equation*}
\sup_{\substack{v \in H^1_0(\Omega)}} \frac{\langle \trace v , \mu_h \rangle}{\|v\|_{H^1(\Omega)}}
\leq C_T \sup_{\substack{v \in H^1_0(\Omega)}} \frac{\langle \trace v , \mu_h \rangle}{ \|\trace v\|_{H^{\frac 12}_{00}(\Gamma)}}.
\end{equation*}
Using Lemma \ref{lemma:prob1_orthproj} we obtain,
\begin{multline*}
C_T \sup_{\substack{v \in H^1_0(\Omega)}} \frac{\langle \trace v , \mu_h \rangle}{ \|\trace v\|_{H^{\frac 12}_{00}(\Gamma)}}= C_T \sup_{\substack{v \in H^1_0(\Omega)}} \frac{\langle P_h(\trace v) , \mu_h \rangle}{ \|\trace v\|_{H^{\frac 12}_{00}(\Gamma)}}
\\
\leq  C \sup_{\substack{v \in H^1_0(\Omega)}} \frac{\langle P_h(\trace v) , \mu_h \rangle}{ \|P_h(\trace v)\|_{H^{\frac 12}_{00}(\Gamma)}}
= C \sup_{\substack{q_h \in Q_h}} \frac{\langle q_h , \mu_h \rangle}{  \|q_h\|_{H^{\frac 12}_{00}(\Gamma)}}.
\end{multline*}
\end{proof}

\begin{theorem}[Discrete inf-sup] The inequality \eqref{eq:infsup_discrete} holds true, 
namely there exists a positive constant $\beta_{h,1}$ such that,
\begin{equation}\label{inf_sup_discrete_prob1}
\inf_{\mu_h \in Q_h} 
\sup_{\substack{v_h \in X_{h,0}^k(\Omega),\\ {\vd}_h \in X_{h,0}^k(\Lambda)}} \frac{ \langle \trace v_h - \ext {\vd}_h, \mu_h \rangle _{\Gamma}} {\vertiii{[v_h, {\vd} _h]} \|\mu_h\|_{H^{-\frac 12 }(\Gamma)}}
\geq \beta_{h,1}. 
\end{equation}
\end{theorem}

\begin{proof}
As in the continuos case, we choose ${\vd}_h =0$ and we have
\begin{equation*}
\sup_{\substack{v_h \in X_{h,0}^k(\Omega),\\ {\vd}_h \in X_{h,0}^k(\Lambda)}} \frac{ \langle \trace v_h - \ext {\vd}_h, \mu_h \rangle _{\Gamma}} {\vertiii{[v_h,  {\vd}_h]}}
\geq \sup_{v_h \in X_{h,0}^k(\Omega)} \frac{ \langle \trace v_h, \mu_h \rangle _{\Gamma} } {\|v_h\|_{H^1(\Omega)}}.
\end{equation*}
Using Lemma \ref{lemma:trspace_infsup} and the boundedness of the harmonic extension operator $\mathcal{H}_\Omega$ from $H^{\frac 12}_{00}(\Gamma)$ to $H^1_0(\Omega)$ introduced in the previous section, we have
\begin{equation*}
\gamma \|\mu_h\|_{H^{-\frac 12}(\Gamma)} \leq  \sup_{q_h \in Q_h} \frac{ \langle q_h, \mu_h \rangle _{\Gamma} } {\|q_h\|_{H^{\frac 12}_{00}(\Gamma)}} 
\leq 
\|\mathcal{H}_\Omega\| \sup_{q_h \in Q_h} \frac{ \langle q_h, \mu_h \rangle _{\Gamma} } {\|\mathcal{H}_\Omega q_h\|_{H^1(\Omega)}} .
\end{equation*}
Let $R_h: H^1_0(\Omega) \rightarrow X_{h,0}^k(\Omega)$ be a quasi interpolation operator (such as the Scott-Zhang operator) satisfying 
$\|R_h v\|_{H^1(\Omega)} \leq C_R \|v\|_{H^1(\Omega)}$ for any  $v \in H^1_0(\Omega)$.
Therefore, we obtain 
\begin{equation*}
\|\mathcal{H}_\Omega\| \sup_{q_h \in Q_h} \frac{ \langle q_h, \mu_h \rangle _{\Gamma} } {\|\mathcal{H}_\Omega q_h\|_{H^1(\Omega)}} 
\leq
\|\mathcal{H}_{\Omega}\| C_R \sup_{q_h \in Q_h} \frac{ \langle q_h, \mu_h \rangle _{\Gamma} } {\|R_h \mathcal{H}_{\Omega} q_h\|_{H^1(\Omega)}}\,.
\end{equation*}
Now we use the conformity of $\mathcal{T}^{\Omega}_h$ to the interface $\Gamma$
to guarantee that the operator $\trace R_h  \mathcal{H}_{\Omega}$ coincides with the identity on the space $Q_h$. Then, for any $q_h \in Q_h$ we have $q_h = \trace R_h  \mathcal{H}_{\Omega} q_h$
and owing to this property we obtain the following inequality, which proves the condition, with $\beta_{h,1} = \gamma \|\mathcal{H}_{\Omega}\|^{-1} C_R^{-1}$,
\begin{multline*}
\gamma \|\mu_h\|_{H^{-\frac 12}(\Gamma)} 
\leq 
\sup_{q_h \in Q_h} \frac{ \langle q_h, \mu_h \rangle_{\Gamma} } {\|q_h\|_{H^{\frac 12}_{00}(\Gamma)}} 
\leq
\|\mathcal{H}_{\Omega}\| C_R \sup_{q_h \in Q_h} \frac{ \langle q_h, \mu_h \rangle_{\Gamma} } {\|R_h \mathcal{H}_{\Omega} q_h\|_{H^1(\Gamma)}}
\\
=
\|\mathcal{H}_{\Omega}\| C_R \sup_{q_h \in Q_h} \frac{ \langle \trace R_h  \mathcal{H}_{\Omega}q_h, \mu_h \rangle_{\Gamma} } {\|R_h \mathcal{H}_{\Omega} q_h\|_{H^1(\Omega)}} 
\leq \|\mathcal{H}_{\Omega}\| C_R \sup_{v_h \in X_{h,k}(\Omega)} \frac{ \langle \trace v_h, \mu_h \rangle_{\Gamma} } {\|v_h\|_{H^1(\Omega)}}\,. 
\end{multline*}
\end{proof}


\subsubsection{Problem 3D-1D-1D}
In this case, we use the same spaces $X_{h,0}^k(\Omega)$, $X_{h,0}^k(\Lambda)$ defined previously.
For the multiplier space we choose $Q_h=X_{h,0}^k(\Lambda)$, therefore we impose homogeneous Dirichlet boundary condition on $\Lambda \cap \partial \Omega$ also for the Lagrange multiplier. 
We aim to find  $u_h \in X_{h,0}^k(\Omega) ,\ {\ud}_h \in X_{h,0}^k(\Lambda), \ \ldh \in Q_h \subset H^{-\frac12}(\Lambda)$, such that
\begin{subequations}
\begin{align}\label{eq:prob2_discrete}
&(u_h,v_h)_{H^1(\Omega)} + ({\ud}_h,{\vd}_h)_{H^1(\Lambda), |{\cal D}|} 
+  \langle  \mtrace v_h -  {\vd}_h, \ldh \rangle_{\Lambda, |\DD|} 
\\
\nonumber
&\quad = (f,v_h)_{L^2(\Omega)} + (\avrd{g},{\vd}_h)_{L^2(\Lambda),|\D|}
\quad \forall v_h \in X_h(\Omega), \ {\vd}_h \in X_h(\Lambda)\,,
\\
& \langle \mtrace u_h - {\ud}_h, \mdh \rangle_{\Lambda,| \DD| } 
= \paolo{\langle \avrc{q},\mdh \rangle_{\Lambda,|\DD|}}
\quad \forall \mdh \in Q_h\,.
\end{align}
\end{subequations}
Below we address the well-posedness of the 3D-1D-1D discrete problem with this alternative choice of multiplier space.

\begin{lemma}\label{lemma:prob2_orthproj}
Let $P_h: H^{\frac 12}_{00}(\Lambda) \longrightarrow Q_h$ be the orthogonal projection operator defined  for any $v \in H^{\frac 12}_{00}(\Lambda)$ by $(P_h v , \psi)_{\Lambda,|\DD|}= (v, \psi)_{\Lambda , |\DD|} \qquad \forall \psi \in Q_h$.  
Then, $P_h$ is continuous on $H^{\frac 12}_{00}(\Lambda)$, namely
$\|P_h v\|_{H^{\frac 12}_{00}(\Lambda),|\DD|} \leq C \|v\|_{H^{\frac 12}_{00}(\Lambda),|\DD|}$.
\end{lemma}

\begin{lemma}\label{infsup_avr_trspace}
There exist a constant $\gamma >0$ such that for any $\mu_h \in Q_h$,
\begin{equation*}
\sup_{\substack{q_h \in Q_h}} \frac{\langle q_h , \mdh \rangle_{\Lambda, |\DD|}}{ \|q_h\|_{H^{\frac 12}_{00}(\Lambda),|\DD|}} \geq \gamma \|\mdh\|_{H^{-\frac 12}(\Lambda)}\,.
\end{equation*} 
\end{lemma}

The proofs of these Lemmas follow the ones of Lemmas \ref{lemma:prob1_orthproj} and \ref{lemma:trspace_infsup} with the only difference that the arguments are applied to $\Lambda$ instead of $\Gamma$.

\begin{theorem}[Discrete inf-sup] 
The inequality \eqref{eq:infsup_discrete} holds, 
namely there exists a positive constant $\beta_{h,2}$ such that,
\begin{equation}
\inf_{\mu_h \in Q_h} 
\sup_{\substack{v_h \in X_{h,0}^k(\Omega),\\ {\vd}_h \in X_{h,0}^k(\Lambda)} }\frac{\langle \mtrace v_h -  {\vd}_h, \mdh \rangle _{\Lambda,|\DD|} } {\vertiii{[v_h, {\vd}_h]} \|\mdh\|_{H^{-\frac 12 }(\Lambda)} } 
\geq \beta_{h,2}. 
\end{equation}
\end{theorem}
\begin{proof}
 Again, we choose ${\vd}_h =0$, so that the proof reduces to showing that there exists $\beta_{h,2}$ such that
\begin{equation*}
\sup_{v_h \in X_{h,0}^k(\Omega)} \frac{ \langle \mtrace v_h , \mdh \rangle _{\Lambda,|\DD|} } {\|v_h\|_{H^1(\Omega)} }\geq \beta_{h,2} \|\mdh\|_{H^{-\frac 12}(\Lambda)} 
\quad \forall \mdh \in Q_h.
\end{equation*}
For any $w \in H^{\frac 12}(\Lambda)$, Lemma \ref{lemma:H12norm} ensures that 
$\|\ext w\|_{H^{\frac 12}_{00}(\Gamma)} = \|w\|_{H^{\frac 12}_{00}(\Lambda),|\DD|}$.
As in the previous case, we use the extension operator $\mathcal{H}_{\Omega}$ from $H^{\frac 12}_{00}(\Gamma)$ to $H^1_0(\Omega)$ and the quasi interpolation operator $R_h$ from $H^1_0(\Omega)$ to $X_{h,0}^k(\Omega)$. Then, we exploit the conformity of the meshes on $\Omega$, $\Gamma$ and $\Lambda$ and the fact that the operator $\mtrace R_h \mathcal{H}_{\Omega} \ext$ coincides with the identity if applied to functions in $Q_h$. As a result, from Lemma \ref{infsup_avr_trspace}, we obtain the following inequality 
\begin{align*}
\gamma \|\mdh\|_{H^{-\frac 12}(\Lambda)} &\leq 
\sup_{q_h \in Q_h} \frac{ \langle q_h, \mdh \rangle_{\Lambda,|\DD|} } {\|q_h\|_{H^{\frac 12}_{00}(\Lambda),|\DD|}} 
=  \sup_{q_h \in Q_h} \frac{ \langle q_h, \mdh \rangle _{\Lambda,|\DD|}} {\|\ext q_h\|_{H^{\frac 12}_{00}(\Gamma)}} 
\\
&\leq \|\mathcal{H}_{\Omega}\| \sup_{q_h \in Q_h} \frac{ \langle q_h, \mdh \rangle _{\Lambda,|\DD|} } {\|\mathcal{H}_{\Omega} \ext q_h\|_{H^1(\Omega)}} 
\leq C_R\|\mathcal{H}_{\Omega}\| \sup_{q_h \in Q_h} \frac{ \langle q_h, \mdh \rangle _{\Lambda,|\DD|} } {\|R_h \mathcal{H}_{\Omega} \ext q_h\|_{H^1(\Omega)}}
\\
&=  C_R\|\mathcal{H}_{\Omega}\|\sup_{q_h \in Q_h} \frac{ \langle \mtrace R_h \mathcal{H}_{\Omega} \ext q_h, \mdh \rangle _{\Lambda,|\DD|}} {\|R_h \mathcal{H}_{\Omega} \ext w_h\|_{H^1(\Omega)}} \\
&\leq C_R\|\mathcal{H}_{\Omega}\| \sup_{v_h \in X_{h,0}^k} \frac{ \langle \mtrace v_h, \mdh \rangle _{\Lambda,|\DD|}} {\|v_h\|_{H^1(\Omega)}}\,,
\end{align*}
that concludes the proof with $\beta_{h,2} = \gamma \|\mathcal{H}_{\Omega}\|^{-1} C_R^{-1}$.
\end{proof}

\def\patch{\omega _j}
\subsection{Analysis of the case where $\mathcal{T}^{\Omega}_h$ and $\mathcal{T}^{\Lambda}_\fkh$ do not conform to $\Gamma$} \label{sec:unfit2}
We analyze now the case in which the elements of the 3D mesh $\mathcal{T}^{\Omega}_h$
do not conform with the surface $\Gamma$ nor with $\Lambda$. 
As the 3D-1D-1D formulation is more suitable for this purpose,
we solely focus on the analysis of the discrete version of Problem 3D-1D-1D.

\subsubsection{Problem 3D-1D-1D} 
Let $u_h \in X_{h,0}^1(\Omega)$ be the approximation of the 3D problem
and let $\udfh \in X_{\fkh,0}^1(\Lambda)$ the one of the 1D problem.
In contrast to the conforming case, 
here we limit the analysis to the case of piecewise-linear finite elements. 
With little abuse of notation, we use the sub-index $h$ for the product space 
$X_h = X_{h,0}^1(\Omega) \times X_{\fkh,0}^1(\Lambda)$.
Concerning the multiplier space, let $\mathcal{G}_h = \{K \in \mathcal{T}^{\Omega}_h: \, K\cap \Lambda \neq \emptyset\}$, be the set of the 3D elements that intersect $\Lambda$. Then we define $Q_h=\{\ldh : \ldh \in P^0(K) \, \forall K \in \mathcal{G}_h\}$. We notice that the multiplier functions are defined on the 3D elements. Again with a little abuse of notation, we denote with $Q_h$ also the restriction to $\Lambda$ of the space of piecewise constant functions defined in 3D. As a result, we have $Q_h \subset L^2(\Lambda) \subset H^{-\frac12}(\Lambda)$.
However, with this choice of multipliers the problem is not inf-sup stable, therefore the idea is to add a stabilization term $s(\ldh, \mdh): \ Q_h \times Q_h \rightarrow \mathbb{R}$ to \eqref{eq:prob2_discrete} following the approach introduced in \cite{burman2014}. 
The objective of this section is to analyze the stabilized version of the 3D-1D-1D problem: find $[u_h, \kent{\udfh]} \in X_h$ and $\ldh \in Q_h$ such that
\begin{multline}\label{eq:prob2_stabilized}
a([u_h, \udfh], [v_h, \vdfh]) +
b([v_h, \vdfh], \ldh)
+ b([u_h, \udfh], \mdh)
\\
- s_h(\ldh, \mdh) 
= c(v_h) + d(\mdh)
  \quad \forall [v_h, \vdfh] \in X_h, \,
  \forall \mdh \in Q_h\,.
\end{multline}
The idea of the stabilization strategy proposed in \cite{burman2014}
is to identify a new multiplier space $Q_H$, which is never implemented in practice, such that inf-sup stability with $X_h$ holds true. Then, the stabilization operator is designed to control the distance between $Q_h$ and $Q_H$ through the following inequality
\begin{equation*}
    \|\mdh - \pi_H \mdh\|_{Q_H} \leq C s_h(\mdh, \mdh)\,,
\end{equation*}
being $\pi_H$ a suitable projection operator $Q_h \rightarrow Q_H$. Applying the results obtained in \cite{burman2014}, the well posedness of problem \eqref{eq:prob2_stabilized} is governed by the following lemma.
\begin{lemma}[Lemma 2.3 of \cite{burman2014}]\label{lemma23burman}
\begin{enumerate}
    \item If the $b(\cdot,\cdot),\, X_h,Q_H \rightarrow \R$ is inf-sup stable.
    \item If the stabilization operator $s_h(\cdot,\cdot),\, Q_h,Q_h \rightarrow \R$ is such that
    \begin{equation*}
        \beta \|\mdh\|_{H^{-\frac12}(\Lambda)} \leq \sup\limits_{v_h \in X_h} \frac{b(v_h,\mdh)}{\vertiii{v_h}} + s_h(\mdh, \mdh), \quad \forall \mdh \in Q_h\,.
    \end{equation*}
    \item If for any $[v_h, \vdfh] \in X_h$ there exists a function $\xi_h\in Q_h$ depending on $[v_h, \vdfh]$, namely $\xi_h=\xi_h([v_h, \vdfh])$, s.t.
\begin{gather}
\label{stab_coercivity}
a([v_h, \vdfh],[v_h, \vdfh] ) + b([v_h, \vdfh], \xi_h) \geq \alpha_\xi \vertiii{[v_h, \vdfh]}_{X_h}\,,
\\
\label{stab_stability}
(s_h(\xi_{h}, \xi_{h}))^{\frac 12} \leq c_s \vertiii{[v_h, \vdfh]}_{X_h}\,,
\end{gather}
being $\vertiii{[\cdot \, ,\, \cdot] }_{X_h}$ a suitable discrete norm. 
\end{enumerate}
Then, problem \eqref{eq:prob2_stabilized} admits a unique solution.
\end{lemma}
For the proof of this result we refer the reader to Lemma 2.3 of \cite{burman2014}. In the remainder of this section, we show how to find a multiplier space $Q_H$ and a stabilization operator $s_h$ such that all the assumptions of Lemma \ref{lemma23burman} are satisfied.

\begin{figure}
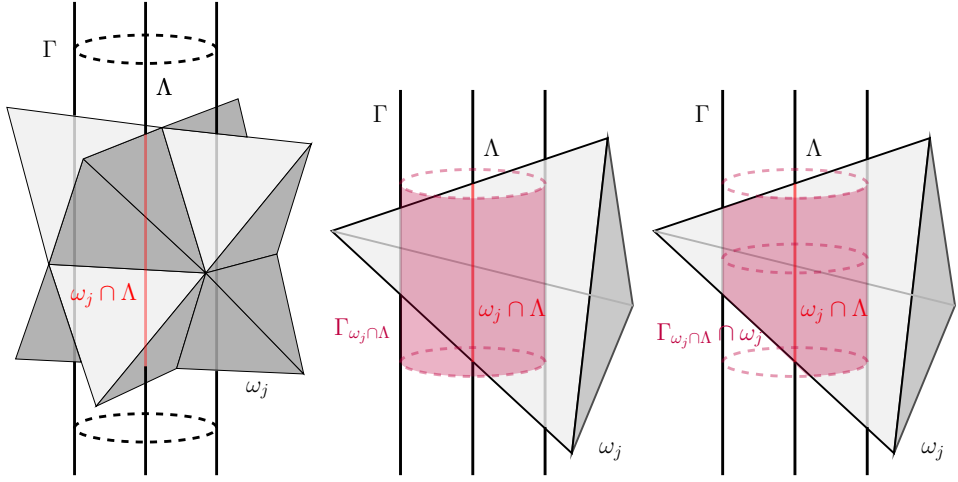
\label{fig:gamma_generated}
\centering
\includegraphics[width=0.32\textwidth]{ext_patch.tex}
\includegraphics[width=0.32\textwidth]{gamma_generated.tex}
\includegraphics[width=0.32\textwidth]{intersection.tex}
\caption{
(Left) Extended patches $\omega_j$. (Middle)
$\Gamma_{\omega_j \cap \Lambda}$, the portion of $\Gamma$ generated by $\omega_j \cap \Lambda$. (Right) the intersection between $\Gamma_{\omega_j \cap \Lambda}$ and $\omega_j$. Here for simplicity $\omega_j$ is represented as a single tetrahedron but actually it is a collection of tetrahedra as shown in left panel.}
\end{figure}

The first step consists of showing
that there exists a discrete space $Q_H$ that satisfies the first assumption of Lemma \ref{lemma23burman}. 
We recall that in the case of Problem 3D-1D-1D, 
\begin{equation*}
b([u_h, \vdfh], \mdh) = \left(\mtrace v_h - \vdfh, \mdh\right)_{\Lambda,|\DD|}.
\end{equation*}
The construction of the inf-sup stable space $Q_H$ is based on macro elements of diameter $H$, where $H$ is sufficiently large. In particular, we assume that there exists positive constants $c_h$ and $c_H$ such that $c_h h\leq H \leq c_H^{-1}h$. The space is constructed assembling the 3D elements of $\mathcal{G}_h$ into macro patches $ \patch $ such that $H \leq |\patch\cap \Lambda|\leq cH$ with $H=\min_{j}|\patch\cap \Lambda|$ and $c\geq 1$. Let $M_j$ be the number of elements of the patch $\omega_j$, namely, $\patch=\cup _{i=0}^{M_j} K_i$, where $K_i \in \mathcal{G}_h$. We assume that $M_j$ is uniformly bounded in $j$ by some $M\in \mathbb{N}$ and that the interiors of the patches $\patch$ are  disjoint.
We define $Q_H$ as the space of piecewise-constant functions on the patches, namely
$Q_H=\left\{{\md }_H: \, {\md }_H\in P^0(\patch)\, \forall j \right\}$.
As previously pointed out for $Q_h$, we denote with $Q_H$ also the restriction of the multiplier space to $\Lambda$, namely say $Q_H \subset L^2(\Lambda) \subset H^{-\frac12}(\Lambda)$.
Moreover, we associate to each patch $\patch$ a shape-regular extended patch (using the classical definition of shape-regularity, see for example \cite{MR2050138}), still denoted by $\omega_j$ for notational simplicity, which is built adding to $\patch$ a sufficient number of elements of $\mathcal{T}_h^{\Omega}$ and we assume that the interiors of the new extended patches $\omega _j$ are still disjoint (see Figure \ref{fig:gamma_generated}). The extended patches $\omega _j$ are built such that they fulfill the conditions meas$(\patch)=\mathcal{O}(H^3)$ and diam$(\Gamma_{\patch\cap \Lambda}\cap \omega_j)=\mathcal{O}(H)$ ($\mathcal{O}(X)$ means $cX \leq \mathcal{O}(X) \leq CX$), where $\Gamma_{\patch\cap \Lambda}$ is the portion of $\Gamma$ with centerline $\patch\cap \Lambda$. The latter assumption is required to ensure that the intersection of $\Gamma_{\patch\cap \Lambda}$ and $\omega_j$ is not too small and it will be needed later on to prove the inf-sup stability of the space $Q_H$ in Lemma \ref{lemma:Lh_infsup}. A representation of this construction in the simple case in which $\omega _j$ is composed just by one tetrahedron is shown in Figure \ref{fig:gamma_generated}.
Thanks to the shape regularity of these extended patches, the following discrete trace inequality holds true for any function $v\in H^1(\omega_j)$, 
\paolo{
\begin{equation}\label{discr_trace_ineq}
\|\trace v\|_{L^2(\Gamma\cap \omega_j)} \leq C_I H^{-\frac 12} \|v\|_{L^2(\omega_j)}
\end{equation}
}
Moreover, $\forall u_h \in X^1_{h,0}(\Omega)$ we have the following average inequality, 
which is a consequence of the definition of $\mtrace$, Jensen inequality, \kent{and the fact that the patches are disjoint}
\begin{multline}\label{eq:avrg_ineq}
\sum _j \|\mtrace u_h\|^2_{L^2(\patch \cap \Lambda),|\DD| }
= \int_{\Lambda} |\DD| \left( \frac{1}{|\DD|} \int_{\DD} \trace u_h \right)^2
\\
\leq \int_{\Lambda} \int_{\DD} (\trace u_h)^2
= \int_{\Gamma} (\trace u_h)^2
= \sum _j \int_{\patch \cap \Gamma} (\trace u_h)^2
= \sum _j \|\trace u_h\|^2_{L^2(\patch \cap \Gamma)}.
\end{multline}



We are now ready to prove that the space $Q_H$ is inf-sup stable. 
\begin{lemma}\label{lemma:Lh_infsup}
The space $Q_H$ is inf-sup stable, 
namely there exists $\beta >0$ such that
\begin{equation*}
\sup_{\substack{v_h \in X_{h,0}^1(\Omega),\\ \vdfh \in X_{\fkh,0}^1(\Lambda)}} \frac{\left(\mtrace v_h - \vdfh, {\md} _H\right)_{\Lambda, |\DD|}}{\vertiii{[v_h, \vdfh]}} \geq \beta \|{\md}_H\|_{H^{-\frac 12}(\Lambda)} \quad \forall {\md}_H \in Q_H.
\end{equation*} 
\end{lemma}

\begin{proof}
We choose $\vdfh=0$ and we prove that
\begin{equation*} 
\sup_{v_h \in X_{h,0}^1(\Omega)} \frac{\left(\mtrace v_h ,{\md}_H\right)_{\Lambda, |\DD|}}{\|v_h\|_{H^1(\Omega)}} \geq \beta \|{\md}_H\|_{H^{-\frac 12}(\Lambda)}.
\end{equation*} 
Proving the last inequality is equivalent to finding the Fortin operator $\pi_F: H^1_0(\Omega) \rightarrow X_{h,0}^1(\Omega)$, such that 
\begin{gather}
\left(\mtrace v - \mtrace \pi _F v  , {\md}_H\right)_{\Lambda, |\DD|}=0, \quad \forall v\in H^1_0(\Omega), \, {\md}_H \in Q_H\,,
\\
\label{eq:cont_Fortin}
\|\pi_F v\|_{H^1(\Omega)} \leq C \|v\|_{H^1(\Omega)}\,.
\end{gather} 

We define
\begin{equation*}
\pi_F v = I_h v + \sum _j \alpha _j \varphi _j \qquad \text{with }\alpha_j =\frac{\int_{\patch \cap \Lambda}|\DD| (\mtrace v-\mtrace I_h v)}{\int_{\patch \cap \Lambda}|\DD|\mtrace \varphi _j}
\end{equation*}
where $I_h: H^1(\Omega) \rightarrow X_{h,0}^1(\Omega)$ denotes an $H^1(\Omega)$-stable interpolant
and $\varphi_j \in X_{h,0}^1(\Omega)$ is such that supp$(\varphi_j)\subset \omega_j$, supp$(\trace \varphi _j) \subset \Gamma_{\patch\cap \Lambda}\cap \omega_j$, $\varphi_j =0$ on $\partial \omega _j$ and 
\begin{equation}\label{phi_properties}
\int_{\patch\cap \Lambda}|\DD|\mtrace \varphi_j=\mathcal{O}(H) \text{ and } \|\nabla \varphi _j\|_{L^2(\omega _j)}=\mathcal{O}(1). 
\end{equation}
We notice that supp$(\trace \varphi_j) \subset \Gamma_{\patch\cap \Lambda}\cap \omega_j$ ensures that $\mtrace \varphi_j \subset \omega_j\cap \Lambda$. Therefore, since the interiors of $\omega_j\cap \Lambda$  are disjoint and $\varphi _j = 0 $ on $\partial \omega_j$,  the functions $\mtrace \varphi_j\, \forall j$ have all disjoint supports. Provided $H$ is sufficiently larger that $h$, the functions $\varphi _j$ and their traces $\trace \varphi _j$ have a sufficiently large support thanks to the fact that $\mathrm{meas}(\patch) = \mathcal{O}(H^3)$ and $\mathrm{diam}(\Gamma_{\patch\cap \Lambda} \cap \omega_j)=\mathcal{O}(H)$.
Owing to these properties it is possible to satisfy \eqref{phi_properties}. 
\kent{Then, by construction,}
\begin{multline*}
\left(\mtrace v - \mtrace \pi _F v  , {\md}_H\right)_{\Lambda, |\DD|} = \sum _j \int_{\patch\cap \Lambda} |\DD |\left[ \mtrace v-\mtrace I_h v-\sum _i \alpha_i \mtrace \varphi _i \right]{\md}_H \\
=\sum _j \int_{\patch\cap \Lambda}|\DD| \left[ \mtrace v-\mtrace I_h v-\alpha_j \mtrace \varphi _j \right]{\md}_H
\\
=\sum _j \left[\int_{\patch\cap \Lambda} |\DD| (\mtrace v-\mtrace  I_h v) {\md}_H
-[\int_{\patch\cap \Lambda} |\DD| (\mtrace v-\mtrace  I_h v) {\md}_H \right]=0.
\end{multline*}
Concerning the continuity of $\pi_F$, we exploit the assumptions that the interiors of $\patch$ are disjoint, $\mathrm{supp}(\varphi _j) \subset \patch$ and the $H^1$-stability of $I_h$ to show that
\begin{equation*}
\|\nabla \pi_F v \|_{L^2(\Omega)} 
\leq C \|\nabla  v\|_{L^2(\Omega)} + \left(\sum_j\alpha_j^2\|\nabla \varphi _j\|^2_{L^2(\omega_j)}\right)^{\frac 12}\,.
\end{equation*}
For the second term,
using that $\|\nabla \varphi _j\|_{{L^2(\omega_j)}}=\mathcal{O}(1)$, 
$\int_{\patch\cap \Lambda}|\DD|\mtrace \varphi_j=\mathcal{O}(H)$ and that
$|\patch\cap \Lambda| \leq cH$,
\kent{exploiting Jensen's average inequality} \eqref{eq:avrg_ineq} 
and trace inequality \eqref{discr_trace_ineq},
and finally applying the approximation properties of $I_h$,
the following upper bound holds true
\paolo{(where all the constants have been condensed into $C$),}
\begin{multline*}
\sum_j \alpha_j ^2\|\nabla \varphi _j\|^2_{L^2(\omega_j)}
\leq C \sum_j \frac{\left(\int_{\patch\cap \Lambda} |\DD| (\mtrace v-\mtrace I_h v)\right)^2}{\left(\int_{\patch\cap \Lambda}|\DD|\mtrace \varphi_j \right)^2}
\\
\leq \frac {C}{H^2} \sum_j |\patch\cap \Lambda| \int_{\patch\cap \Lambda} |\DD|^2(\mtrace v-\mtrace I_h v)^2
\\
\leq \frac {C}{H} \sum_j \| \mtrace (v-I_h v)\|^2_{L^2(\patch\cap \Lambda), |\DD|}
\leq \frac {C}{H} \sum_j \| \trace(v-I_h v)\|^2_{L^2(\omega _j\cap \Gamma)}  
\\
\leq \frac {C}{H^2} \sum_j  \| v-I_h v\|^2_{L^2(\omega_j)} 
\leq C \frac {1}{H^2}  \| v-I_h v\|^2_{L^2(\Omega)} 
\leq C \|\nabla  v\|^2_{L^2(\Omega)}
\end{multline*}
that is the $H^1$-stability of $\pi_F$. We notice that the constant in the inequality \eqref{eq:cont_Fortin} is independent of how $\Lambda$  cuts the elements of the mesh $\mathcal{T}_h^{\Omega}$.
\end{proof}

For the second assumption of Lemma \ref{lemma23burman}, we recall that $b(v_h,\mdh)$ is continuous \kent{with respect to the norms} $\vertiii{v_h},\,\|\mdh\|_{L^2(\Lambda)}$. 
Using Lemma \ref{lemma:Lh_infsup}, and in particular the existence of a Fortin projector,
there exists a constant $\beta$ such that 
(the proof is \miro{analogous} to the one of Lemma 2.1 in \cite{burman2014})
\begin{equation}\label{relax_infsup}
\beta \|\mdh\|_{H^{-\frac12}(\Lambda)} \leq \sup\limits_{v_h \in X_h} \frac{b(v_h,\mdh)}{\vertiii{v_h}} + \|\mdh-\pi_H \mdh\|_{L^2(\Lambda)}, \quad \forall \mdh \in Q_h\,.
\end{equation}
We define $\pi_H = \sum_j \pi_H^j\,: L^2(\Lambda) \rightarrow Q_H$, 
where $\pi_H^j$ is the operator
\begin{equation}\label{eq:projector}
\pi_H^j w_{|\patch \cap \Lambda} =\frac{1}{|\Gamma_{\patch\cap \Lambda}|}\int_{\patch\cap \Lambda}|\DD| w\quad \forall j.
\end{equation}   
Since $\cup_j \omega_j \cap \Lambda  = \Lambda$ and $\omega_j \cap \Lambda$ are not overlapping,
we obtain that $\pi_H$ is an orthogonal projection, namely $(w-\pi_H w,\pi_H w)=0$.
Moreover, for any $w \in L^2(\Lambda)$ the following Poincarè inequality holds true, see for example \cite[Corollary B.65]{MR2050138},
\begin{equation}\label{disc_poincare_ineq}
\|w - \pi_H w\|_{L^2(\omega_j\cap\Lambda),|\DD|} \leq C_P H \|\partial_s w\|_{L^2(\omega_j\cap\Lambda),|\DD|}\,.
\end{equation}
We consider the following stabilization operator
\begin{equation}\label{eq:stab}
s_h(\ldh, \mdh)= \sum _{K\in \mathcal{G}_{h}} \int_{\partial K\setminus \partial \mathcal{G}_{h}} h \llbracket \ldh \rrbracket \llbracket \mdh \rrbracket,
\end{equation}
being $\llbracket \ldh \rrbracket$ the jump of $\ldh$ across the internal faces of $\mathcal{G}_h$.
Then, we use the result of \cite{burman2014}, Section III to show that
\begin{equation*}
    \|\mdh-\pi_H \mdh\|_{L^2(\Lambda)} \leq C s_h(\mdh, \mdh)\,,
\end{equation*}
which combined with \eqref{relax_infsup} shows that the second assumption of Lemma \ref{lemma23burman} holds true.

The third step of the analysis consists of showing that \eqref{stab_coercivity} and \eqref{stab_stability} are satisfied. 
We introduce the following discrete norms
\begin{equation*}
\|\, \lambda \,\|_{\pm \frac 12, h, \Lambda} = \|h^{\mp\frac 12} \lambda\|_{L^2(\Lambda)},
\end{equation*}
recalling that $h$ is the mesh size of $\mathcal{T}^\Omega_{h}$. We equip the space $X_h$ with the discrete norm
\begin{equation*}
\vertiii{[u_h, \udfh]}^2_{X_h}
= \|u_h\|^2_{H^1(\Omega)}+\|\udfh\|^2_{H^1(\Lambda),|\D|} + \|\mtrace u_h - \udfh\|^2_{\frac 12, h, \Lambda, |\DD|},
\end{equation*}
and the space $Q_H$ with the $L^2$ norm $\|{\md}_H \|_{L^2(\Lambda)}$.

Also, the function $\xi_h([v_h, \vdfh]) \in Q_H \subset Q_h \subset L^2(\Lambda)$ is defined as follows
\begin{equation*}
{\xi_h}_{|\patch\cap \Lambda}=\frac{\delta}{H} \pi_H(\mtrace u_h-\udfh)_{|\patch \cap \Lambda}.
\end{equation*}
where $\delta$ is an arbitrarily small parameter. Then the following result holds true. 
\begin{lemma}
Given $\pi_H,\ s_h(\cdot,\cdot), \ \xi_h$ defined above, 
choosing $\delta$ small enough,
the inequalities \eqref{stab_coercivity} and \eqref{stab_stability} are satisfied.
\end{lemma}
\begin{proof} 
Concerning the coercivity property \eqref{stab_coercivity}, we show that $\forall [u_h, \udfh]$, there exists $\xi_h \in Q_h$ such that,
\begin{equation*}
(u_h,u_h)_{H^1(\Omega)}+ (\udfh, \udfh)_{H^1(\Lambda), |\D|} +   (\mtrace u_h - \udfh, \xi_h)_{\Lambda,|\DD|} \geq \alpha_{\xi}\vertiii{[u_h, \udfh]}^2_{X_h}.
\end{equation*}
Using the definitions of $\pi_H$ and $\xi_h([u_h, \udfh])$ previously presented and recalling that $\xi_h\in Q_H \subset Q_h$, we obtain
\begin{multline*}
\left( \mtrace u_h - \udfh, \xi _h \right)_{\Lambda,|\DD|} 
= \frac{\delta}{H}\sum_j \pi_H^j( \mtrace u_h - \udfh)\int_{\patch\cap \Lambda} |\DD|(\mtrace u_h - \udfh)
\\
= \frac{\delta}{H} \sum_j \int_{\patch\cap \Lambda}|\DD| (\pi_H( \mtrace u_h - \udfh) )^2
= \frac{\delta}{H}\sum_j \|\pi_H( \mtrace u_h - \udfh)\|_{L^2(\patch\cap \Lambda),|\DD|}^2
\\
= \frac{\delta}{H} \sum_j \left( \|\mtrace u_h - \udfh\|^2_{L^2(\patch\cap \Lambda),|\DD|} - \|(\pi_H-\mathcal{I})(\mtrace u_h - \udfh)\|^2_{L^2(\patch\cap \Lambda),|\DD|} \right) 
\\ 
\geq \frac{\delta}{H} \sum_j \left(
\|\mtrace u_h-\udfh\|^2_{L^2(\patch\cap \Lambda),|\DD|}
- \|(\pi_H - \mathcal{I})\mtrace u_h\|^2_{L^2(\patch\cap \Lambda), |\DD|} \right.
\\
\left.- \|(\pi_H - \mathcal{I})\udfh\|^2_{L^2(\patch\cap \Lambda),|\DD|} \right). 
\end{multline*}
Now, we seek for an upper bound of the second and third (negative) terms of the last inequality.
For the second term, we apply the additional assumption that the operators $\mtrace$ and $\partial_s$ commute. This is true if the cross section $\D$ does not depend on the arclength $s$.
\paolo{
Then, we use the Poincar\'e inequality \eqref{disc_poincare_ineq}, the average inequality \eqref{eq:avrg_ineq} and the trace inequality \eqref{discr_trace_ineq} to show that,
\begin{multline*}
\sum _j  \|(\pi_H - \mathcal{I})\mtrace u_h\|^2_{L^2(\patch\cap \Lambda),|\DD|} 
\leq C_P^2 H^2 \sum _j \|\mtrace \partial_s  u_h\|_{L^2(\patch\cap \Lambda),|\DD|}^2
\\
\leq C_P^2 H^2 \sum _j \|\trace \partial_s u_h\|_{L^2(\omega _j \cap \Gamma)}^2
\leq C_P^2 C_I^2 H \sum _j \|\nabla u_h\|_{L^2(\omega _j)}^2.
\end{multline*}
For the third term, the following upper bound holds true,
\begin{multline*}
\sum _j \|(\pi_H - \mathcal{I})\udfh\|^2_{L^2(\patch\cap \Lambda),|\DD|} \leq C_P^2 H^2 \sum _j \|\partial_s \udfh\|_{L^2(\patch\cap \Lambda),|\DD|}^2
\\
\leq C_P^2 H^2 \frac{\max |\DD|}{\min |\D|} \sum _j \|\partial_s \udfh\|_{L^2(\patch\cap \Lambda),|\D|}^2.
\end{multline*}
Combining the last three inequalities, reminding that $c_h h\leq H \leq c_H^{-1}h$, 
we obtain
\begin{multline*}
a([u_h, \udfh],[u_h, \udfh] ) + b([u_h, \udfh], \xi_h([u_h, \udfh]))
\geq (1-\delta C_P^2 C_I^2) \|\nabla u_h\|^2_{L^2(\Omega)} 
\\
+ \left(1- \delta C_P^2 H \frac{\max |\DD|}{\min |\D|}\right) \|\partial_s \udfh\|^2_{L^2(\Lambda), |\D|}
+\delta c_H  \|\mtrace u_h-\udfh\|^2_{\frac 12,h,\Lambda, |\DD|}
\end{multline*}
and choosing $\delta=\frac{1}{2}\min\left[(C_P^2 C_I^2)^{-1}, \left(C_P^2 H \frac{\max |\DD|}{\min |\D|}\right)^{-1}\right]$ we obtain the desired inequality.
Concerning inequality \eqref{stab_stability}, the proof is analogous to the one in \cite{burman2014}.}
\end{proof}


\section{A benchmark problem with analytical solution}

Let $\Omega=[0,1]^3$, $\Lambda=\{x=\tfrac{1}{2}\}\times \{y=\tfrac{1}{2}\} \times [0,1] $
and $\Omega_{\ominus}=[\tfrac{1}{4}, \tfrac{3}{4}]\times [\tfrac{1}{4}, \tfrac{3}{4}]\times [0, 1]$.
As a benchmark for the two formulations we consider the following coupled problems
\begin{subequations}\label{benchmark}
\begin{align}
\label{benchm_3d}
-\Delta u=f \quad &\text{in $\Omega$},\\
\label{benchm_1d}
-d_{zz}^2 \ud =\avrd{g} \quad &\text{on $\Lambda$},\\
u=u_b \quad &\text{on $\partial \Omega$},
\end{align}
\end{subequations}
where for formulation \eqref{eq:problem1} the mix-dimensional coupling constraint reads
\begin{equation}
  \label{eq:couple_1}
\trace{u} - \ext{\ud} = q_1\quad\text{ on }\Gamma,
\end{equation}
while for \eqref{eq:problem2} we set
\begin{equation}
    \label{eq:couple_2}
\avrc{u} - \ud = \avrc{q}_2\quad\text{ on }\Lambda.
\end{equation}
In \eqref{benchmark}-\eqref{eq:couple_2} the right-hand sides shall be defined as 
\begin{eqnarray*}
  &f=8\pi ^2 \sin (2\pi x) \sin (2\pi y),\quad &\avrd{g}={\pi ^2}\sin \left({\pi z}\right),\quad u_b=\sin (2\pi x) \sin (2\pi y),\\
  &q_1=\sin (2\pi x) \sin (2\pi y) - \sin \left({\pi z}\right),\quad &\avrc{q}_2=-\sin \left({\pi z}\right).
\end{eqnarray*}

The exact solution of \eqref{benchmark}, regardless of the coupling constraint,
is given by
\begin{equation}
\label{benchm_sol3d-1d}
u=\sin (2\pi x) \sin (2\pi y),
\quad
\ud=\sin \left({\pi z}\right).
\end{equation}
Let us notice that $\ud$ satisfies homogeneous Dirichlet conditions at the boundary of $\Lambda$.
Moreover, the solution \eqref{benchm_sol3d-1d} satisfies on $\Gamma$ the relation
\begin{equation}\label{benchm_flux}
\nabla u \cdot \textbf{n}_{\oplus}=d_z \ud n_{\oplus,z}=0,
\end{equation}
with $n_{\oplus,z}$ the $z-$component of the normal unit vector to $\Gamma$.

\paolo{We prove that the solution of \eqref{benchmark} is equivalent to the one of \eqref{eq:problem2}.}
Precisely, we prove that \eqref{benchm_sol3d-1d} is the solution of
\eqref{eq:problem2}. Using the integration by part formula and homogeneous
boundary conditions on $\Omega$ and $\Lambda$, from \eqref{eq:problem2} we have
\begin{align*}
&-(\Delta u, v)_{L^2(\Omega)} - |{\cal D}|(d^2_{ss} \ud, \vd)_{L^2(\Lambda)} 
+ |{\cal D}|\langle \avrc{v}  - \vd, \ld \rangle_\Lambda
\\
\nonumber
&\qquad\qquad= (f,v)_{L^2(\Omega)} + |{\cal D}| (\avrd{g},\vd)_{L^2(\Lambda)}
\quad \forall v \in H^1_0(\Omega), \vd \in H^1(\Lambda).
\end{align*}
Since $\ld=0$ and the first of \eqref{benchm_sol3d-1d} satisfies \eqref{benchm_3d} and the second satisfies \eqref{benchm_1d}, we have that
\begin{align*}
-(\Delta u, v)_{L^2(\Omega)} =  (f,v)_{L^2(\Omega)}, 
\\
-|{\partial \cal D}|(d^2_{ss} \ud, \vd)_{L^2(\Lambda)}  = |{\cal D}| (\avrd{g},\vd)_{L^2(\Lambda)}.
\end{align*}
Thus \eqref{benchm_sol3d-1d} satisfy equations \eqref{eq:problem2a}, \eqref{eq:problem2b}.
The fact that the solution satisfy \eqref{eq:problem2c} follows from \eqref{eq:couple_2}.

We can prove in a similar way that \eqref{benchm_sol3d-1d}, with $\lambda=0$
satisfy \eqref{eq:problem1}. Note in particular that $q_1$ is such that
$\trace{u} - \ext{\ud} = q_1$ on $\Gamma$.

\subsection{Numerical experiments. $\mathcal{T}^{\Omega}_h$ conforming to $\Gamma$}\label{sec:experiment_conform}
\begin{figure}
 \centering
 \includegraphics[width=0.45\textwidth]{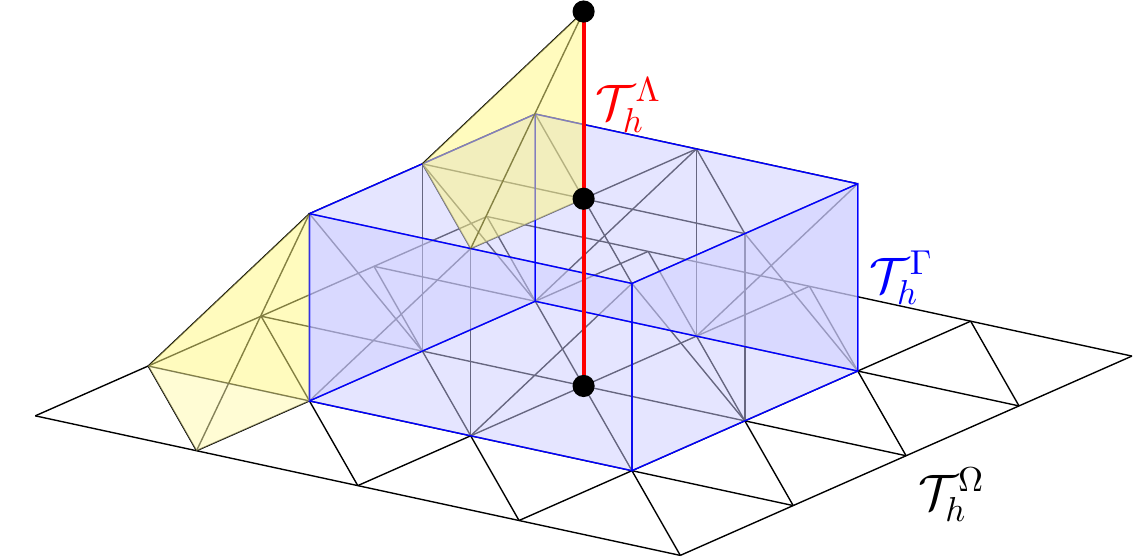}
 \includegraphics[width=0.45\textwidth]{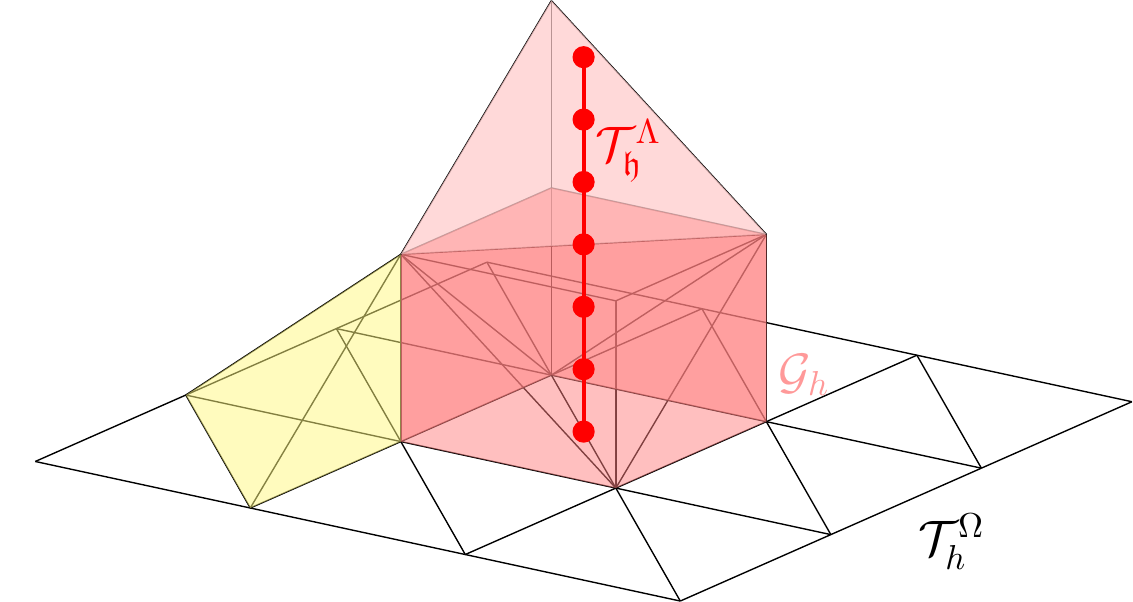} 
 \vspace{-10pt}
 \captionof{figure}{
   (Left)
$\Lambda$ and $\Gamma$ conforming discretization of $\Omega$
   used for \eqref{eq:problem1} and \eqref{eq:problem2}.
   (Right) Sample discretization of the benchmark geometry in the non-conforming case for \eqref{eq:problem2}. 
 }
\label{fig:mesh}
\end{figure}

Using the benchmark problem \eqref{benchmark} we now investigate convergence
properties of the two formulations. To this end we consider a \emph{uniform} mesh
of $\mathcal{T}^{\Omega}_h$ of $\Omega$ consisting of tetrahedra with diameter $h$.
Further, the discretization shall be geometrically \emph{conforming} to both $\Lambda$
and $\Gamma$ such that the meshes $\mathcal{T}^{\Gamma}_h$, $\mathcal{T}^{\Lambda}_h$
are made up of facets and edges of $\mathcal{T}^{\Omega}_h$ respectively, cf. Figure \ref{fig:mesh}
for illustration.

\begin{table}
  \scriptsize{
  \begin{center}
    \begin{tabular}{l|llll}
\hline
\multicolumn{5}{c}{\textcolor{blue}{$\mathcal{T}^{\Omega}_h$ conforming to $\Gamma$, $\Lambda$}}\\
\hline
    $h^{-1}$ & $\norm{u-u_h}_{H^1(\Omega)}$ & $\norm{\ud-\udh}_{H^1(\Lambda)}$ & $\norm{\lambda-\lambda_h}_{H^{-1/2}(\Gamma)}$ & $\norm{\lambda-\lambda_h}_{L^2(\Gamma)}$\\
      \hline
4  & 3.4E0(--)    & 5.3E-1(--)   & 2.9E0(--)    &8.7E0(--)    \\
8  & 1.7E0(0.99)  & 2.6E-1(1.06) & 6.1E-1(2.25) &1.9E0(2.21)  \\
16 & 8.7E-1(0.99) & 1.3E-1(1.02) & 1.4E-1(2.13) &4.7E-1(1.99) \\
32 & 4.4E-1(1.00) & 6.3E-2(1.00) & 3.4E-2(2.03) &1.3E-1(1.80) \\
64 & 2.2E-1(1.00) & 3.1E-2(1.00) & 8.6E-3(2.00) &4.2E-2(1.68) \\
\hline
$h^{-1}$ & $\norm{u-u_h}_{H^1(\Omega)}$ & $\norm{\ud-u_{\odot}}_{H^1(\Lambda)}$ & $\norm{\ld-\ldh}_{H^{-1/2}(\Lambda)}$ & $\norm{\ld-\ldh}_{L^2(\Lambda)}$\\
\hline
4   & 3.1E0(--)    & 5.4E-1(--)   & 4.4E-2(--)   & 7.8E-2(--)  \\
8   & 1.7E0(0.87)  & 2.6E-1(1.06) & 1.1E-2(2.01) & 1.9E-2(2.01)\\
16  & 8.6E-1(0.96) & 1.3E-1(1.02) & 2.7E-3(2.01) & 4.8E-3(2.02)\\
32  & 4.4E-1(0.99) & 6.3E-2(1.00) & 6.7E-4(2.01) & 1.2E-3(2.01)\\
64  & 2.2E-1(1.00) & 3.1E-2(1.00) & 1.7E-4(2.01) & 3.0E-4(2.01)\\
128 & 1.1E-1(1.00) & 1.6E-2(1.00) & 4.1E-5(2.01) & 7.4E-5(2.00)\\
\hline
\multicolumn{5}{c}{\textcolor{red}{$\mathcal{T}^{\Omega}_h$ non conforming to $\Gamma$, $\Lambda$}}\\
\hline
    $h^{-1}$ & $\norm{u-u_h}_{H^1(\Omega)}$ & $\norm{\ud-u_{\odot\mathfrak{h}}}_{H^1(\Lambda)}$ & \multicolumn{2}{c}{$\norm{\ld-\ldh}_{L^2(\mathcal{G}_h)}$}\\
      \hline
5   & 2.6E0(--)    & 2.3E-1(--)   & \multicolumn{2}{c}{1.7E-1(--)} \\ 
9   & 1.5E0(0.84)  & 9.4E-2(1.42) & \multicolumn{2}{c}{7.1E-2(1.36)}\\
17  & 8.1E-1(0.94) & 4.3E-2(1.18) & \multicolumn{2}{c}{2.9E-2(1.37)}\\
33  & 4.2E-1(0.98) & 2.1E-2(1.06) & \multicolumn{2}{c}{7.9E-3(1.91)}\\
65  & 2.1E-1(0.99) & 1.1E-2(1.02) & \multicolumn{2}{c}{2.6E-3(1.64)}\\
129 & 1.1E-1(1.00) & 5.2E-3(1.01) & \multicolumn{2}{c}{8.5E-4(1.61)}\\
\hline
    \end{tabular}
  \end{center}    
  }
  \captionof{table}{Error convergence on a benchmark problem \eqref{benchmark}.
    (Top) problem \eqref{eq:problem1}, (middle) \eqref{eq:problem2} with
    conforming discretization and (bottom) \eqref{eq:problem2} in case
    $\mathcal{T}^{\Omega}_h$ does not conform to $\Lambda$ using 
    stabilized formulation \eqref{eq:prob2_stabilized}. Continuous linear Lagrange
    elements are used for $u_h$, $\udh$ and $u_{\odot \mathfrak{h}}$ and $\ldh$ in
    conforming case, while in nonconforming case $\ldh$ is piecewise constant
    on elements of $\mathcal{G}_h$.
  }
  \label{tab:error_conform}
\end{table}

Considering inf-sup stable discretization in terms of continuous linear Lagrange
($P_1$) elements (for all the spaces), Table \ref{tab:error_conform}
lists the errors of formulations \eqref{eq:problem1} and \eqref{eq:problem2}
on the benchmark problem. It can be seen that the error in $u$ and $\ud$ in $H^1$ norm
converges linearly (as can be expected due to $P_1$ element discretization).
Moreover, the error of the Lagrange multiplier approximation in $H^{-1/2}$ norm
decreases quadratically. In the light of $P_1$ discretization this rate appears
superconvergent. We speculate that the result is due to the fact that the
exact solution is particularly simple, $\lambda=\ld=0$.
We remark that for $u$ and $\ud$ the error is interpolated into the finite element space of
piecewise quadratic \emph{discontinous} functions. For \eqref{eq:problem2} we
evaluate the fractional norm and interpolate the error using piecewise continuous
cubic functions. 
For the sake of comparison with non-conforming formulation of \eqref{eq:problem2} from
\S\ref{sec:unfit2} Table \ref{tab:error_conform} also
lists the error of the Lagrange multiplier in the $L^2$ norm. Here, quadratic convergence is observed
for \eqref{eq:problem2}. For \eqref{eq:problem1} the rate is between 1.5 and 2.

\begin{figure}
\centering
\includegraphics[width = 0.6\textwidth]{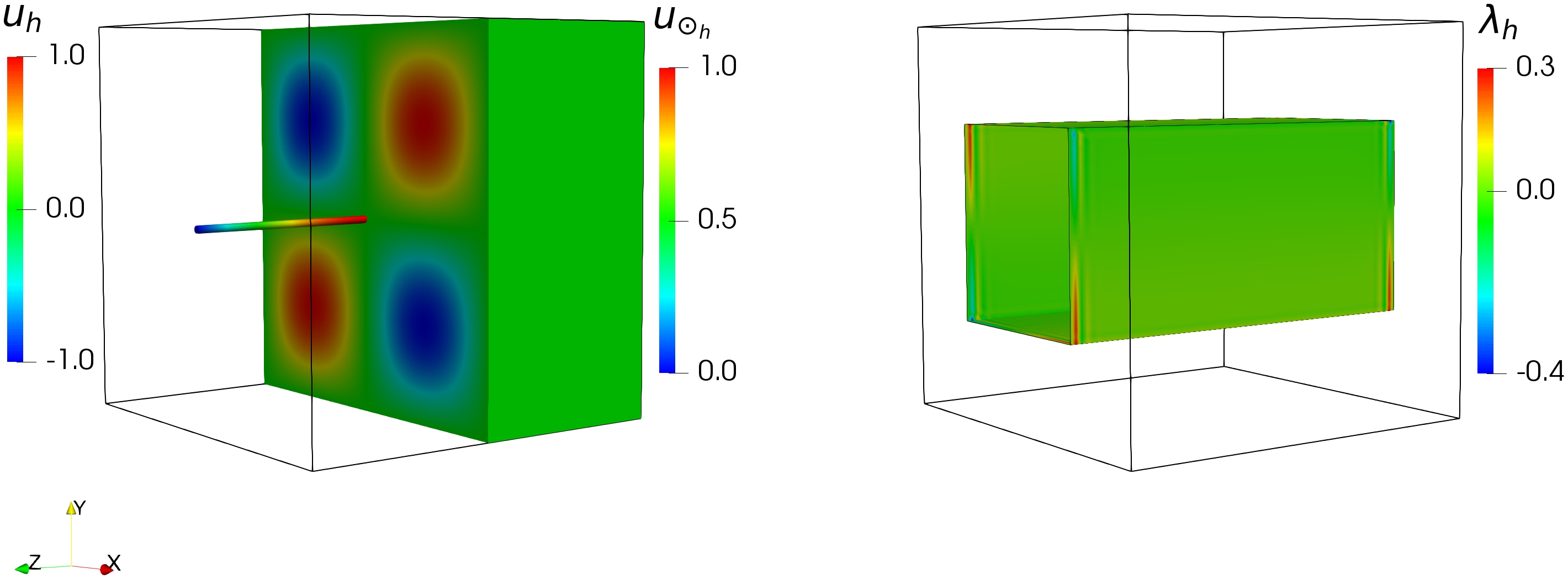}
\includegraphics[width = 0.3\textwidth]{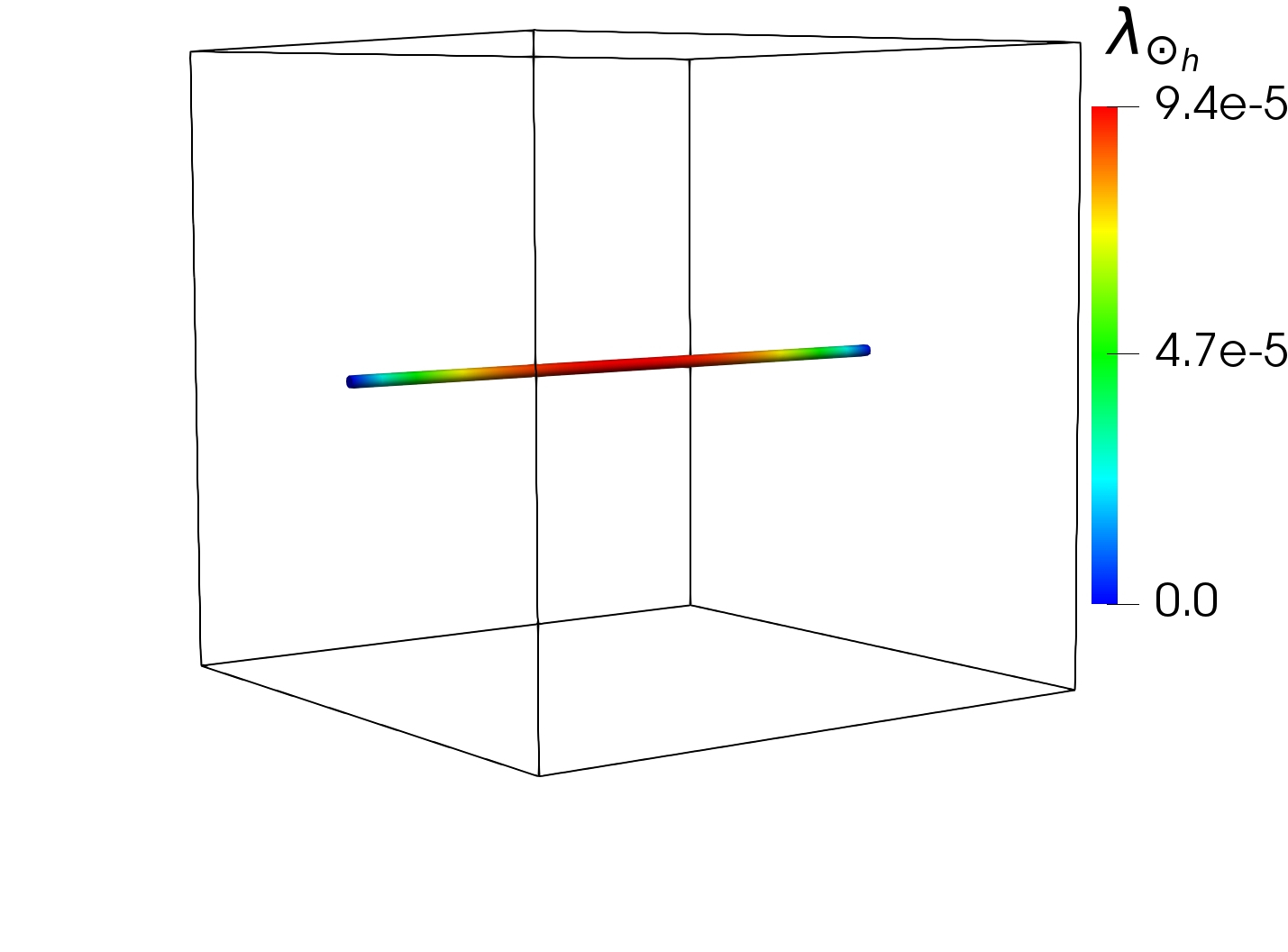}
\vspace{-10pt}
\caption{Numerical solution of problem \eqref{eq:problem1} and
  \eqref{eq:problem2}. (Left) functions $u_h$, $\udh$ (practically identical
  in both problems). (Middle) Lagrange multiplier for \eqref{eq:problem1}
  and (right) for \eqref{eq:problem2}.}
\label{fig:sol_benchm1}
\end{figure}

We plot the numerical solution of problem \eqref{eq:problem1} and \eqref{eq:problem2} in
Figure \ref{fig:sol_benchm1}.

\subsection{Numerical experiments. $\mathcal{T}^{\Omega}_h$ non-conforming to $\Gamma$}\label{sec:experiment_nonconform}
Using benchmark problem \eqref{benchmark} we consider \eqref{eq:problem2} in the
setting of \S \ref{sec:unfit2}. To this end we let $\mathcal{T}^{\Omega}_h$ be a uniform
mesh of $\Omega$ such that no cell $\mathcal{T}^{\Omega}_h$ has any edge
lying on $\Lambda$. Further we let $\mathfrak{h}=h/3$ in $\mathcal{T}^{\Lambda}_{\mathfrak{h}}$,
cf. Figure \ref{fig:mesh}.

Using discretization in terms of $P_1$-$P_1$-$P_0$ element Table \ref{tab:error_conform}
lists the error of the stabilized formulation of \eqref{eq:problem2}. A linear
convergence in the $H^1$ norm can be observed in the error of $u$ and $\ud$. We
remark that the norms were computed as in \S\ref{sec:experiment_conform}. For simplicity
the convergence of the multiplier is measured in the $L^2$ norm rather then the $H^{-1/2}(\Gamma)$
norm used in the analysis. Then, convergence exceeding order 1.5 can be observed, however,
the rates are rather unstable.


\subsection{Comparison}
In Tables \ref{tab:error_conform} one can observe that 
all the formulations yield practically identically accurate approximations of $u$.
Further, compared to the conforming case, the stabilized formulation \eqref{eq:problem2}
results in a greater accuracy of $u_{\odot h}$ as the underlying mesh $\mathcal{T}^{\Lambda}_{h}$ is
here finer. Due to the different definitions in the three formulations, comparision of the Lagrange
multiplier convergence is not straightforward. We therefore limit ourselves to a
comment that in the $L^2$ norm all the formulations yield faster than linear convergence.
In order to discuss solution cost of the formulations we consider 
the resulting preconditioned linear systems. In particular, we shall compare
spectral condition numbers and the time to convergence of the preconditioned
minimal residual (MinRes) solver with the with stopping criterion requiring
the relative preconditioned residual norm to be less than $10^{-8}$. We remark
that we shall ignore the setup cost of the preconditioner.
Following operator preconditioning technique \cite{mardal2011preconditioning} we
propose as preconditioners for \eqref{eq:problem1} and \eqref{eq:problem2} in the
conforming case the (approximate) Riesz mapping with respect to the inner products of
the spaces in which the two formulations were proved to be well posed.
In particular, the preconditioner for the Lagrange multiplier relies on
(the inverse of) the fractional Laplacian $-\Delta^{-1/2}$ on $\Gamma$ for
\eqref{eq:problem1} and $\Lambda$ for \eqref{eq:problem2}.
A detailed analysis of the preconditioners will be presented in a separate
work. We remark that in both cases the fractional Laplacian was here realized
by spectral decomposition \cite{kuchta2016preconditioners}.
For the unfitted stabilized formulation \eqref{eq:problem2} the Lagrange multiplier preconditioner
uses a Riesz map with respect to the inner product due to $L^2(\mathcal{G}_h)$ and
the stabilization \eqref{eq:stab}, i.e.
\[
(\ldh, \mdh) \mapsto \sum _{K\in \mathcal{G}_h}\int_{K}\ldh \mdh + \sum _{K\in \mathcal{G}_h} \int_{\partial K \setminus \partial \mathcal{G}_h} h \llbracket \ldh \rrbracket \llbracket \mdh \rrbracket.
\]
This simple choice does not yield bounded
iterations. However, establishing a robust preconditioner in this case 
is beyond the scope of the paper and shall be pursued in the future works.
In Table \ref{tab:cost} we compare solution time, number of iterations and
condition numbers of the (linear systems due to the) three formulations.
Let us first note that the proposed preconditioners for \eqref{eq:problem1} and
\eqref{eq:problem2} in the conforming case seem robust with respect to discretization
parameter as the iteration counts and condition numbers are bounded in $h$.
We then see that the solution time for \eqref{eq:problem1} is about 2 times
longer compared to \eqref{eq:problem2} which is about 4 times more expensive
than the solution of the Poisson problem \eqref{benchm_3d}. This is in addition to the
higher setup costs of the preconditioner, which in our implementation involve solving
an eigenvalue problem for the fractional Laplacian. Therefore it is advantageous to keep
the multiplier space as small as possible. We remark that the missing
results for \eqref{eq:problem1} in Table \ref{tab:cost} are due to the
memory limitations encountered when solving the eigenvalue problem
for the Laplacian, which for finest mesh involves cca 32 thousand eigenvalues, cf. Appendix \ref{sec:appendix}.
Due to the missing proper preconditioner for the Lagrange multiplier block the
number of iterations in the third, unfitted formulation can be seen to approximately
double on refinement. 

\begin{table}
  \scriptsize{
  \begin{center}
    \begin{tabular}{l|lll|lll|lll|ll}
      \hline
      \multirow{2}{*}{$l$} & \multicolumn{3}{c|}{\eqref{eq:problem1}} & \multicolumn{3}{c|}{\eqref{eq:problem2}} & \multicolumn{3}{c|}{ Stabilized \eqref{eq:problem2}} & \multicolumn{2}{c}{\eqref{benchm_3d}}\\
      \cline{2-12}
      & \# & $T\left[s\right]$ & $\kappa$ & \# & $T\left[s\right]$ & $\kappa$ & \# & $T\left[s\right]$ & $\kappa$ & \# & $T\left[s\right]$\\
      \hline
      1 & 20 & 0.03  & 15.56 & 9  & 0.02  & 3.04 & 21  & 0.01   & 9.70  &3 & $<0.01$\\ 
      2 & 35 & 0.06  & 16.28 & 17 & 0.03  & 4.67 & 31  & 0.03   & 15.87 &4 & $<0.01$\\ 
      3 & 38 & 0.14  & 16.64 & 22 & 0.06  & 6.25 & 53  & 0.15   & 32.93 &5 & 0.01   \\ 
      4 & 39 & 1.70  & 16.75 & 24 & 0.89  & 7.03 & 110 & 4.54   & 61.48 &5 & 0.12   \\ 
      5 & 38 & 12.04 & 16.78 & 20 & 5.21  & 5.02 & 232 & 59.43  & 94.25 &5 & 0.90  \\ 
      6 & -- & --    & --    & 17 & 28.77 & --   & 507 & 832.90 & --    &6 & 7.75  \\
      \hline
    \end{tabular}
    \end{center}
    }
  \caption{Cost comparison of the formulations across refinement levels $l$.
    Number of Krylov  iterations (preconditioned conjugate gradient for \eqref{benchm_3d},
    MinRes otherwise) and the condition number of the preconditioned
    problem is denoted by $\#$ and $\kappa$ respectively. Time till convergence
    of the iterative solver (excluding the setup) is shown as $T$. 
  }
\label{tab:cost}
\end{table}

\bibliographystyle{siamplain}
\bibliography{3d1d_coupled}

\begin{thebibliography}{10}

\bibitem{adams1975pure}
{\sc R.~A. Adams}, {\em Pure and applied mathematics 65}, Sobolev Spaces,
  (1975).

\bibitem{MR359352}
{\sc I.~Babu\v{s}ka}, {\em The finite element method with {L}agrangian
  multipliers}, Numer. Math., 20 (1972/73), pp.~179--192.

\bibitem{berg2020modelling}
{\sc M.~Berg, Y.~Davit, M.~Quintard, and S.~Lorthois}, {\em Modelling solute
  transport in the brain microcirculation: is it really well mixed inside the
  blood vessels?}, Journal of Fluid Mechanics, 884 (2020).

\bibitem{Bertoluzza201897}
{\sc S.~Bertoluzza, A.~Decoene, L.~Lacouture, and S.~Martin}, {\em Local error
  estimates of the finite element method for an elliptic problem with a {D}irac
  source term}, Numerical Methods for Partial Differential Equations, 34
  (2018), pp.~97--120.

\bibitem{Nordbotten}
{\sc W.~Boon, J.~Nordbotten, and J.~Vatne}, {\em Mixed-dimensional elliptic
  partial differential equations}, tech. report, arXiv, Cornell University
  Library, 2017.
\newblock arXiv:1710.00556v2.

\bibitem{MR365287}
{\sc F.~Brezzi}, {\em On the existence, uniqueness and approximation of
  saddle-point problems arising from {L}agrangian multipliers}, Rev.
  Fran\c{c}aise Automat. Informat. Recherche Op\'{e}rationnelle S\'{e}r. Rouge,
  8 (1974), pp.~129--151.

\bibitem{burman2014}
{\sc E.~Burman}, {\em Projection stabilization of {L}agrange multipliers for
  the imposition of constraints on interfaces and boundaries}, Numerical
  Methods for Partial Differential Equations, 30 (2014), pp.~567--592.

\bibitem{Cerroni2019}
{\sc D.~Cerroni, F.~Laurino, and P.~Zunino}, {\em Mathematical analysis, finite
  element approximation and numerical solvers for the interaction of 3d
  reservoirs with 1d wells}, GEM - International Journal on Geomathematics, 10
  (2019).

\bibitem{DAngelo}
{\sc C.~D'Angelo}, {\em {Multi scale modelling of metabolism and transport
  phenomena in living tissues, PhD Thesis}}, EPFL, Lausanne, 2007.

\bibitem{d2012finite}
{\sc C.~D'Angelo}, {\em Finite element approximation of elliptic problems with
  {D}irac measure terms in weighted spaces: applications to one-and
  three-dimensional coupled problems}, SIAM Journal on Numerical Analysis, 50
  (2012), pp.~194--215.

\bibitem{d2008coupling}
{\sc C.~D'Angelo and A.~Quarteroni}, {\em On the coupling of 1d and 3d
  diffusion-reaction equations: Application to tissue perfusion problems},
  Mathematical Models and Methods in Applied Sciences, 18 (2008),
  pp.~1481--1504.

\bibitem{MR2050138}
{\sc A.~Ern and J.-L. Guermond}, {\em Theory and practice of finite elements},
  vol.~159 of Applied Mathematical Sciences, Springer-Verlag, New York, 2004.

\bibitem{fang2008oxygen}
{\sc Q.~Fang, S.~Sakad{\v{z}}i{\'c}, L.~Ruvinskaya, A.~Devor, A.~M. Dale, and
  D.~A. Boas}, {\em Oxygen advection and diffusion in a three dimensional
  vascular anatomical network}, Optics express, 16 (2008), p.~17530.

\bibitem{MR1940257}
{\sc I.~Gansca, W.~F. Bronsvoort, G.~Coman, and L.~Tambulea}, {\em
  Self-intersection avoidance and integral properties of generalized
  cylinders}, Comput. Aided Geom. Design, 19 (2002), pp.~695--707.

\bibitem{gjerde2019singularity}
{\sc I.~G. Gjerde, K.~Kumar, and J.~M. Nordbotten}, {\em A singularity removal
  method for coupled 1d--3d flow models}, Computational Geosciences,  (2019),
  pp.~1--15.

\bibitem{gould2017capillary}
{\sc I.~G. Gould, P.~Tsai, D.~Kleinfeld, and A.~Linninger}, {\em The capillary
  bed offers the largest hemodynamic resistance to the cortical blood supply},
  Journal of Cerebral Blood Flow \& Metabolism, 37 (2017), pp.~52--68.

\bibitem{koppl2016local}
{\sc T.~K{\"o}ppl, E.~Vidotto, and B.~Wohlmuth}, {\em A local error estimate
  for the {P}oisson equation with a line source term}, in Numerical Mathematics
  and Advanced Applications ENUMATH 2015, Springer, 2016, pp.~421--429.

\bibitem{KVWZ}
{\sc T.~K{\"o}ppl, E.~Vidotto, B.~Wohlmuth, and P.~Zunino}, {\em Mathematical
  modeling, analysis and numerical approximation of second-order elliptic
  problems with inclusions}, Mathematical Models and Methods in Applied
  Sciences, 28 (2018), pp.~953--978.

\bibitem{koppl2014optimal}
{\sc T.~K\"oppl and B.~Wohlmuth}, {\em Optimal a priori error estimates for an
  elliptic problem with {D}irac right-hand side}, SIAM Journal on Numerical
  Analysis, 52 (2014), pp.~1753--1769.

\bibitem{KMM2}
{\sc M.~Kuchta, K.-A. Mardal, and M.~Mortensen}, {\em Preconditioning trace
  coupled 3d-1d systems using fractional {L}aplacian}, Numerical Methods for
  Partial Differential Equations, 0.

\bibitem{kuchta2016preconditioners}
{\sc M.~Kuchta, M.~Nordaas, J.~Verschaeve, M.~Mortensen, and K.~Mardal}, {\em
  Preconditioners for saddle point systems with trace constraints coupling 2d
  and 1d domains}, SIAM Journal on Scientific Computing, 38 (2016),
  pp.~B962--B987.

\bibitem{laurino_m2an}
{\sc {Laurino, F.} and {Zunino, P.}}, {\em Derivation and analysis of coupled
  pdes on manifolds with high dimensionality gap arising from topological model
  reduction}, ESAIM: M2AN, 53 (2019), pp.~2047--2080.

\bibitem{mardal2011preconditioning}
{\sc K.-A. Mardal and R.~Winther}, {\em Preconditioning discretizations of
  systems of partial differential equations}, Numerical Linear Algebra with
  Applications, 18 (2011), pp.~1--40.

\bibitem{Mori2019887}
{\sc Y.~Mori, A.~Rodenberg, and D.~Spirn}, {\em Well-posedness and global
  behavior of the peskin problem of an immersed elastic filament in stokes
  flow}, Communications on Pure and Applied Mathematics, 72 (2019),
  pp.~887--980.

\bibitem{Peaceman1978183}
{\sc D.~Peaceman}, {\em Interpretation of well-block pressures in numerical
  reservoir simulation.}, Soc Pet Eng AIME J, 18 (1978), pp.~183--194.

\bibitem{Peaceman1983531}
{\sc D.~W. Peaceman}, {\em Interpretation of well-block pressures in numerical
  reservoir simulation with nonsquare grid blocks and anisotropic
  permeability.}, Society of Petroleum Engineers journal, 23 (1983),
  pp.~531--543.

\bibitem{Possenti2019101}
{\sc L.~Possenti, G.~Casagrande, S.~Di~Gregorio, P.~Zunino, and M.~Costantino},
  {\em Numerical simulations of the microvascular fluid balance with a
  non-linear model of the lymphatic system}, Microvascular Research, 122
  (2019), pp.~101--110.

\bibitem{sauter1999extension}
{\sc S.~Sauter and R.~Warnke}, {\em Extension operators and approximation on
  domains containing small geometric details}, East West Journal of Numerical
  Mathematics, 7 (1999), pp.~61--77.

\bibitem{secomb2004green}
{\sc T.~W. Secomb, R.~Hsu, E.~Y. Park, and M.~W. Dewhirst}, {\em Green's
  function methods for analysis of oxygen delivery to tissue by microvascular
  networks}, Annals of biomedical engineering, 32 (2004), pp.~1519--1529.

\bibitem{Vinje2019}
{\sc V.~Vinje, G.~Ringstad, E.~Lindstrøm, L.~Valnes, M.~Rognes, P.~Eide, and
  K.-A. Mardal}, {\em Respiratory influence on cerebrospinal fluid flow – a
  computational study based on long-term intracranial pressure measurements},
  Scientific Reports, 9 (2019).

\end{thebibliography}

\newpage

\appendix

\section{Derivation of the model}
This section provides a rigorous derivation of 3D-1D-1D problem \eqref{eq:problem2} and 3D-1D-2D problem \eqref{eq:problem1}. 
The steps are similar to the derivation presented in \cite{laurino_m2an}, 
however, here the coupling conditions are different,
giving rise to coupled problems featuring Lagrange multipliers. 
Precisely, the starting point is the problem arising from \emph{Dirichlet-Neumann} conditions. 
Find $\up,\uf$ s.t.:
\begin{subequations}\label{eq:dirneu}
\begin{align}
\label{eq:dirneu1}
- \Delta \up  + \up &= f  && \text{ in } \Omega_{\oplus},\\
\label{eq:dirneu2}
- \Delta \uf  + \uf &= g  && \text{ in } \Omega_\ominus,\\
\label{eq:dirneu4}
\up - \uf &= q && \text{ on }  \Gamma,\\
\label{eq:dirneu3}
\nabla\left( \up - \uf \right) \cdot \nn_{\oplus} &=0  && \text{ on } \Gamma,\\
\label{eq:dirneu5}
\up &= 0 && \text{ on } \partial \Omega\,.
\end{align}
\end{subequations}
The coupling constraints defined on $\Gamma$ involve essential or strong conditions.
Such conditions will be enforced weakly by using the method of Lagrange multipliers~\cite{MR359352}. Then, the variational formulation of problem \eqref{eq:dirneu}
is to find $\up \in H^1_{\partial\Omega}(\Omega_\oplus), \ \uf \in H^1_{\partial\Omega_\ominus\setminus\Gamma}(\Omega_\ominus), \ \lambda \in H^{-\frac 1 2}(\Gamma)$ s.t.
\begin{subequations}\label{eq:weak_dirneu}
\begin{align}
&(u_\oplus,v_\oplus)_{H^1(\Omega_\oplus)} + (u_\ominus,v_\ominus)_{H^1(\Omega_\ominus)} 
+ \langle  v_\oplus - v_\ominus, \lambda \rangle_{\Gamma} 
\\
\nonumber
&\quad = (f,v_\oplus)_{L^2(\Omega_\oplus)} + (g,v_\ominus)_{L^2(\Omega_\ominus)}
\quad \forall v_\oplus \in H^1_{\partial\Omega}(\Omega_\oplus), \ v_\ominus \in H^1_{\partial\Omega_\ominus\setminus\Gamma}(\Omega_\ominus),
\\
& \langle u_\oplus - u_\ominus, \mu \rangle_{\Gamma} = 0
\quad \forall  \mu \in H^{-\frac 1 2}(\Gamma)\,.
\end{align}
\end{subequations}
where $\lambda$ is the Lagrange multiplier and it is equivalent to $\nabla \uf \cdot \nn_\ominus$.

\subsubsection*{Model reduction of the problem on $\Omega_{\ominus}$}

We apply the averaging technique to equation \eqref{eq:dirneu2}. In particular, we consider an arbitrary portion $\mathcal{P}$ of the cylinder $\Omega_\ominus$, with lateral surface $\Gamma _{\mathcal{P}}$ and bounded by two perpendicular sections to $\Lambda$, namely $\mathcal{D}(s_1), \ \mathcal{D}(s_2)$ with $s_1<s_2$. We have,
\begin{multline*}
\int_{\mathcal{P}} -\Delta \uf + \uf d\omega =
-\int_{\partial \mathcal{P}} \nabla \uf \cdot \nn_{\ominus} \, d\sigma  + \int_{\mathcal{P}}\uf d\omega=
\\
 \int_{\mathcal{D}(s_1)} \partial_s \uf d\sigma -  \int_{\mathcal{D}(s_2)} \partial_s \uf d\sigma -  \int_{\Gamma_{\mathcal{P}}} \nabla \uf \cdot \nn_{\ominus} d\sigma +\int_{\mathcal{P}}\uf d\omega
\end{multline*}
By the fundamental theorem of integral calculus 
\begin{equation*}
\int_{\mathcal{D}(s_1)} \partial_s \uf d\sigma -  \int_{\mathcal{D}(s_2)} \partial_s \uf d\sigma 
= -\int_{s_1}^{s_2} d_s \int_{\mathcal{D}(s)}  \partial_s \uf d\sigma ds
= -\int_{s_1}^{s_2} d_s \left( |\mathcal{D}(s)| \avrd{\partial_s \uf} \right)
\end{equation*}
Moreover, we have
\begin{equation*}
\int_{\Gamma_{\mathcal{P}}} \nabla \uf \cdot \nn_{\ominus} d\sigma =  \int_{\Gamma_{\mathcal{P}}} \lambda \, d\sigma
=  \int_{s_1}^{s_2} \int_{\partial\mathcal{D}(s)} \lambda d\gamma \, ds 
= \int_{s_1}^{s_2} |\DD(s)| \avrc{\lambda} \,ds\,.
\end{equation*}
From the combination of all the above terms with the right hand side, we obtain that the solution $\uf$ of \eqref{eq:dirneu2} satisfies,
\begin{equation*}
\int_{s_1}^{s_2} \left[ 
-d_s (|\D(s)|\avrd{\partial_s \uf}) + |\D(s)|\avrd{u}_{\ominus} 
- |\partial\mathcal{D}(s)| \avrc{\lambda} - |\mathcal{D}(s)| \avrd{g}
\right]\, ds = 0.
\end{equation*}
Since the choice of the points $s_1,s_2$ is arbitrary, we conclude that the following equation holds true,
\begin{equation}\label{eq:averaged_f}
-d_s (|\D(s)|\avrd{\partial_s \uf}) + |\D(s)|\avrd{u}_{\ominus} 
- |\partial\mathcal{D}(s)| \avrc{\lambda}
= |\mathcal{D}(s)| \avrd{g} \quad \text{on} \ \Lambda\,,
\end{equation}
which is complemented by the following conditions at the boundary of $\Lambda$,
\begin{equation}\label{eq:averaged_f_boundary}
|\D(s)| \avrd{\partial_s \uf}= 0, \quad \text{on} \quad s=0,S.
\end{equation}
Then, we consider variational formulation of the averaged equation \eqref{eq:averaged_f}.
After multiplication by a test function $\vd \in H^1(\Lambda)$, integration on $\Lambda$ and suitable application of integration by parts, we obtain,
\begin{multline*}
\int_\Lambda |\D(s)| \avrd{\partial_s \uf} d_s \vd \, ds - (|\D(s)| \avrd{\partial_s \uf}) \vd |_{s=0}^{s=S}- \int_\Lambda |\DD(s)| \avrc{\lambda} \vd \, ds
+\int_\Lambda |\D(s)| \avrd{u}_\ominus \vd 
\\
= \int_\Lambda |\D(s)| \avrd{g} V\, ds\,.
\end{multline*}
Using boundary conditions, we obtain,
\begin{equation}\label{eq:averaged_f_weak}
  (\avrd{\partial_s \uf}, d_s \vd )_{\Lambda,|\D|} 
+ (\avrd{u}_\ominus, \vd)_{\Lambda, |\D|}
- (\avrc{\lambda},\vd)_{\Lambda,|\DD|}
= (\avrd{g},V)_{\Lambda,|\D|}\,.
\end{equation}
Let us now formulate the modelling assumption that allows us to reduce equation \eqref{eq:averaged_f_weak} to a solvable one-dimensional (1D) model.
More precisely, we assume that the function $\uf$ has a \emph{uniform profile} on each cross section $\mathcal{D}(s)$, namely $\uf(r,s,t) = \ud(s)$.
Therefore, observing that $\ud=\avrc{u}_{\ominus}=\avrd{u}_{\ominus}$, 
and that $\avrd{\partial_s \uf} = \avrd{\partial_s \ud} = d_s \ud$,
problem \eqref{eq:averaged_f_weak} turns out to find $\ud \in H^1(\Lambda)$ such that
\begin{equation}\label{eq:1Ddirneu_weak}
 (d_s \ud, d_s \vd)_{\Lambda,|\mathcal{D}|} 
+(\ud, \vd)_{\Lambda, |\D|}
-(\avrc{\lambda},\vd)_{\Lambda, |\DD|}
=  (\avrd{g},\vd)_{\Lambda, |\mathcal{D}|}
\quad \forall \vd \in H^1(\Lambda)\,.
\end{equation}

\subsubsection*{Topological model reduction of the problem on $\Omega_{\oplus}$}
We focus here on the subproblem of \eqref{eq:dirneu1} related to $\Omega_{\oplus}$.
We multiply both sides of \eqref{eq:dirneu1} by a test function $v\in H^1_0(\Omega)$ and integrate on $\Omega_\oplus$. Integrating by parts and using boundary and interface conditions, we obtain
\begin{equation*}
\begin{split}
\int_{\Omega _{\oplus}}fv\,d\omega&=
\int _{\Omega _{\oplus}}\nabla \up \cdot\nabla v\,d\omega -\int _{\partial \Omega _{\oplus}}\nabla \up \cdot \nn_{\oplus} v\,d\sigma + \int_{\Omega_{\oplus}} \up v \,d\omega
\\
&=\int _{\Omega _{\oplus}}\nabla \up \cdot\nabla v\,d\omega -\int_{\Gamma} \nabla \up \cdot \nn_{\oplus} v\,d\sigma + \int_{\Omega_{\oplus}} \up v\,d\omega
\\
&=\int _{\Omega _{\oplus}}\nabla \up \cdot\nabla v\,d\omega +\int_{\Gamma} \lambda v\,d\sigma + \int_{\Omega_{\oplus}} \up v\,d\omega.
\end{split}
\end{equation*} 
Then, we make the following modelling assumption:
we identify the domain $\Omega_{\oplus}$ with the entire $\Omega$, 
and we correspondingly omit the subscript $\oplus$ to the functions defined on $\Omega_{\oplus}$,
namely
\begin{equation*}
\int_{\Omega_{\oplus}} u_\oplus\, d\omega \simeq \int_{\Omega} u\, d\omega\,.
\end{equation*}
Therefore, we obtain 
\begin{equation*}
(\nabla u ,\nabla v)_{\Omega} +(u,v)_{\Omega}+(\lambda, v)_{\Gamma}  =(f,v)_{\Omega}
\end{equation*}
and combining with \eqref{eq:1Ddirneu_weak} we obtain the first formulation of the reduced problem.

Hence, we have obtained the Problem 3D-1D-2D, equation 
\eqref{eq:problem1}:
Find $u \in H^1_0(\Omega), \ \lambda \in H^{-\frac 1 2}(\Gamma ),\ \ud \in H^1_0(\Lambda)$, such that
\begin{subequations}
\begin{align*}
&(u,v)_{H^1(\Omega)} + (\ud,\vd)_{H^1(\Lambda),|\D|} 
+ \langle \trace v  - \ext \vd, \lambda \rangle_\Gamma &&
\\
\nonumber
&\qquad= (f,v)_{L^2(\Omega)} +  (\avrd{g},\vd)_{L^2(\Lambda),|\D|}\,,
\quad &&\forall v \in H^1_0(\Omega), \ \vd \in H^1(\Lambda)\,,
\\
& \langle \trace u - \ext \ud , \mu \rangle_\Gamma 
= \langle q, \mu \rangle_\Gamma\,,
\quad && \forall \mu \in H^{-\frac 1 2}(\Gamma)\,.
\end{align*}
\end{subequations}
This coupled problem is classified as 3D-1D-2D because the unknowns $u$, $\ud$, $\lambda$ belong to $\Omega \subset \R^3$, $\Lambda \subset \R$ and $\Gamma \subset \R^2$ respectively.
Then, we apply a topological model reduction of the interface conditions, namely we go from a 3D-1D-2D formulation
involving sub-problems on $\Omega$ and $\Lambda$ and coupling operators defined on $\Gamma$
to a 3D-1D-1D formulation where the coupling terms are set on $\Lambda$. 
To this purpose, let us write the Lagrange multiplier and the test functions on every cross section $\partial\mathcal{D}(s)$ as their average plus some fluctuation,
\begin{equation*}
\lambda=\avrc{\lambda}+\tilde{\lambda}, \qquad v=\avrc{v}+\tilde{v},
\quad \text{on} \ \partial\mathcal{D}(s)\,,
\end{equation*}
where $\avrc{\tilde{\lambda}}=\avrc{\tilde{v}}=0$. 
Therefore, the coupling term on $\Gamma$ can be decomposed as,
\begin{equation*}
\int_{\Gamma}\lambda v\, d\sigma
=\int _{\Lambda}  \int_{\partial\mathcal{D}(s)} (\avrc{\lambda}+\tilde{\lambda})(\avrc{v}+\tilde{v})d\gamma ds
= \int_{\Lambda}|\partial\mathcal{D}(s)| \avrc{\lambda}\avrc{v}\,ds+\int_{\Lambda}  \int_{\partial\mathcal{D}(s)} \tilde{\lambda}\tilde{v}d\gamma ds\,.
\end{equation*}
Thanks to the additional assumption that the product of fluctuations is small,
\begin{equation*}
\int_{\partial\mathcal{D}(s)} \tilde{\lambda}\tilde{v} d\gamma \simeq 0\,
\end{equation*}
the term $\left(\trace v, \lambda \right)_{\Gamma}$ becomes $\left( \mtrace v, \avrc{\lambda} \right)_{\Lambda, |\DD|}$, where $\mtrace$ denotes the composition of operators $\avrc{(\cdot)}\circ \trace$.  Combined with \eqref{eq:1Ddirneu_weak}, this leads to the 3D-1D-1D formulation of the reduced problem, namely equation \eqref{eq:problem2}:
find $u \in H^1_0(\Omega)$, $\ud \in H^1_0(\Lambda)$, $\ld \in H^{-\frac 1 2}(\Lambda)$, such that
\begin{subequations}
\begin{align*}
&(u,v)_{H^1(\Omega)} + (\ud,\vd)_{H^1(\Lambda),|\D|} 
+  \langle \mtrace v - \vd, \ld \rangle_{\Lambda, |\DD|} &&
\\
\nonumber
&\qquad\qquad= (f,v)_{L^2(\Omega)} +  (\avrd{g},V)_{L^2(\Lambda),|\D|}\,,
\quad && \forall v \in H^1_0(\Omega), \ \vd \in H^1_0(\Lambda)\,,
\\
& \langle \mtrace u - \ud, \md \rangle_{\Lambda,|\DD|} 
= \langle q, \md \rangle_\Gamma = \langle \avrc{q}, \md \rangle_{\Lambda,|\DD|} \,,
\quad &&\forall \md \in H^{-\frac 1 2}(\Lambda)\,.
\end{align*}
\end{subequations}


\section{Proof of Lemma \ref{lemma:H12norm}}

\begin{proof}
Let us consider the eigenvalue problem for the Laplace operator on $\Gamma$ with homogeneous Dirichlet conditions at $x=0,\, X$ and periodic boundary conditions at $y=0,\, Y$. Let us also consider the Laplace eigenproblem on $(0,X)$ with homogeneous Dirichlet conditions. Let us denote as $\phi _{ij}(x,y)$ and $\rho _{ij}$, for $i=1,2,\dots$, $j=0,1,\dots$, the eigenfunctions and the eigenvalues of the Laplacian on $\Gamma$, and with $\phi _i(x)$ and $\rho _i$ the eigenfunctions and the eigenvalues of the Laplacian on $(0,X)$. In particular,
\begin{align*}
&\phi _{ij}(x,y)=\sin \left(\frac{i\pi x}{X}\right)\left( \cos\left(\frac{j2\pi y}{Y}\right)+ \sin\left(\frac{j2\pi y}{Y}\right) \right),\quad
&&\rho_{ij}=\left(\frac{i\pi}{X}\right) ^2+\left(\frac{j2\pi}{Y}\right)^2,
\\
&\phi _{i}(x)=\sin \left(\frac{i\pi x}{X}\right),\quad
&&\rho _i = \left(\frac{i\pi}{X}\right) ^2.
\end{align*}

We use here the following representation of the fractional norms, 
\begin{equation}\label{eq:fracnorm_lambda}
\begin{aligned}
&\|u\|_{H^{\frac 12}_{00}(\Lambda)}=\left(\sum_{i=1}^{\infty} \left( 1+ \rho_{i}\right)^{\frac 12}|a_i|^2\right)^{\frac 12},\\
&\|u\|_{H^{\frac 12}_{00}(\Gamma)}=
\sum_{i=1}^{\infty} \sum_{j=1}^{\infty}  \left( 1+ \left(\frac{i\pi}{X}\right)^2 
+ \left(\frac{j2\pi}{Y}\right)^2\right)^{\frac 12} |a_{i,j}|^2
\end{aligned}
\end{equation} 
with $a_i=\left(u, \phi _i \right)_{\Lambda}$ and $a_{ij}=\left( u, \phi_{ij}\right)_{\Gamma}$.
It is easy to verify that 
\begin{equation}\label{int_eigenf}
\int_0^{Y} \phi _{ij}(x,y)=0 \quad \forall j>0, \forall i\,,
\quad
\int_0^{Y} \phi _{ij}(x,y)= Y \, \sin\left(\frac{i\pi x}{X}\right) \quad \mbox{if } j=0, \forall i  .
\end{equation}
Moreover we recall that $\phi_{i,j}(x,y)$ and $\phi _i(x)$ form an orthogonal basis of $L^2(\Gamma)$ and $L^2(0,X)$ respectively. Therefore,
\begin{equation*}
\avrc{u}(x)=\frac{1}{Y}\int_0^{Y} u(x,y)\, dy
= \frac{1}{Y}\sum_{i,j} a_{i,j}\int_0^{Y}  \phi_{i,j}(x,y) \, dy
=  \sum_{i} a_{i,0} \phi_{i}(x).
\end{equation*}

Let the constant $C$ be equal to 
$C=C(X)=\sum_{i=1}^{\infty}\left( 1+ \left(\frac{i\pi}{X}\right)^2\right)^{\frac 12}$. 
Then,  from \eqref{eq:fracnorm_lambda} we have 
\begin{multline*}
\|\avrc{u}\|^2_{H^{\frac 12}_{00}(0,X)}
=\sum_{i=1}^{\infty} \left( 1+ \rho_{i}\right)^{\frac 12}a_i^2
\\
=C\left( \int_0^X \avrc{u}(x) \sin\left(\frac{i\pi x}{X}\right)\, dx \right)^2
= C\left( \sum_{j=1}^\infty a_{j,0}\int_0^X \sin\left(\frac{j\pi x}{X}\right) \sin\left(\frac{i\pi x}{X}\right) \, dx  \right)
\\
= \sum_{i=1}^{\infty} \frac{X^2}{4} \left( 1+ \left(\frac{i\pi}{X}\right)^2\right)^{\frac 12}a_{i,0}^2
\leq \frac{X^2}{4}\sum_{i=1}^{\infty} \sum_{j=1}^{\infty}  \left( 1+ \left(\frac{i\pi}{X}\right)^2 
+ \left(\frac{j2\pi}{Y}\right)^2\right)^{\frac 12} |a_{i,j}|^2
\\ 
=\frac{X^2}{4}\|u\|^2_{H^{\frac 12}_{00}(\Gamma)},
\end{multline*}
where we have used the orthogonality property
\begin{equation*}
\int_0^X \sin\left(\frac{i\pi x}{X}\right) \sin\left(\frac{j\pi x}{X}\right)\, dx=\begin{cases} 0 & i\neq j\\
\frac X 2 & i = j\\
\end{cases}
\end{equation*}
and we have applied \eqref{eq:fracnorm_lambda} in the last equality.
As a result of the previous inequality, we have proved the first statement of the Corollary, namely $u\in H^\frac12_{00}(\Gamma) \rightarrow \overline{u}\in H^\frac12_{00}(\Lambda)$.

The second statement of the Corollary addresses the case of the function $u$ constant with respect to $y$. Precisely, we have
\begin{equation*}
    \begin{aligned}
&\|u\|^2_{H^{\frac 12}_{00}(\Gamma)}=\sum_{i=1}^{\infty}\sum_{j=0}^{\infty} \left( 1+ \rho_{ij}\right)^{\frac 12}|a_{ij}|^2\\
&=\sum_{i=1}^{\infty}\sum_{j=0}^{\infty} \left(  1+ \left(\frac{i\pi}{X}\right)^2 + \left(\frac{j2\pi}{Y}\right)^2\right)^{\frac 12}\left( \int _0^X\int _0^Y u(x,y )\phi_{ij}(x,y) \right)^2
\\
&=\sum_{i=1}^{\infty}\sum_{j=0}^{\infty} \left(  1+ \left(\frac{i\pi}{X}\right)^2 + \left(\frac{j2\pi}{Y}\right)^2\right)^{\frac 12}\left( \int _0^X u(x) \int _0^Y \phi_{ij}(x,y) \right)^2,
\end{aligned}
\end{equation*}
and using \eqref{int_eigenf} we obtain
\begin{equation*}
\begin{aligned}
\|u\|^2_{H^{\frac 12}_{00}(\Gamma)}
&=\sum_{i=1}^{\infty}\left( 1+ \left(\frac{i\pi}{X}\right)^2\right)^{\frac 12}\left(\int _0^X Yu(x)\, \sin\left(\frac{i\pi x}{X}\right)\right)^2\\
&=Y^2 \sum_{i=1}^{\infty}\left( 1+ \rho _i\right)^{\frac 12}|a_i|^2 = Y^2  \|u\|^2_{H^{\frac 12}_{00}(0,X)}.
\end{aligned}
\end{equation*}
\end{proof}

\section{System sizes in benchmark formulations} \label{sec:appendix}
In Table \ref{tab:dof} we list dimensions of the finite element spaces used to discretize formualations \eqref{eq:problem1}, \eqref{eq:problem2} and stabilized
\eqref{eq:problem2} on different levels of refinement. The number of 
degrees of freedom in subspace $W_{i, h}$ is denote as $\lvert W_{i, h}\rvert$. 
We recall that the discrete spaces  
are $X^1_{h, 0}(\Omega)\times X^1_{h, 0}(\Lambda)\times Q_h(\Gamma)$
for the 3D-1D-2D problem \eqref{eq:problem1},
$X^1_{h, 0}(\Omega)\times X^1_{h, 0}(\Lambda)\times Q_h(\Lambda)$
for the 3D-1D-1D problem \eqref{eq:problem2}, and 
$X^1_{h, 0}(\Omega)\times X^1_{\mathfrak{h}, 0}(\Lambda)\times Q_h(\mathcal{G}_h)$
for the stabilized 3D-1D-1D problem.

\begin{table}[h!]
\begin{center}
  \scriptsize{
    \begin{tabular}{l|lll|lll|lll}
      \hline
      \multirow{2}{*}{$l$} & \multicolumn{3}{c|}{\eqref{eq:problem1}} & \multicolumn{3}{c|}{\eqref{eq:problem2}} & \multicolumn{3}{c}{ Stabilized \eqref{eq:problem2}}\\
      \cline{2-10}
       & $\lvert W_{1, h}\rvert$  & $\lvert W_{2, h}\rvert$ & $\lvert W_{3, h}\rvert$ 
       & $\lvert W_{1, h}\rvert$  & $\lvert W_{2, h}\rvert$ & $\lvert W_{3, h}\rvert$
       & $\lvert W_{1, h}\rvert$  & $\lvert W_{2, \mathfrak{h}}\rvert$ & $\lvert W_{3, h}\rvert$\\
      \hline
     1& 125    & 5  & 40   & 125     & 5   & 5   & 180     & 13  & 24  \\
     2& 729    & 9  & 144  & 729     & 9   & 9   & 900     & 25  & 48  \\
     3& 4913   & 17 & 544  & 4913    & 17  & 17  & 5508    & 49  & 96  \\
     4& 35937  & 33 & 2112 & 35937   & 33  & 33  & 38148   & 97  & 192 \\
     5& 275K & 65 & 8320 & 275K  & 65  & 65  & 283K  & 193 & 384 \\
     6& --     & -- & -    & 2.15M & 129 & 129 & 2.18M & 385 & 768 \\
     \hline
    \end{tabular}
    }
\end{center}
    \caption{Number of degrees of freedom of the discrete spaces used in the numerical experiments.}
    \label{tab:dof}
\end{table}

\end{document}